\documentclass{amsart}
\usepackage{amsthm, amsfonts, amssymb, amscd} 
\newcommand{\C}{\mathbb{C}}

\newcommand{\GL}{\mathrm{GL}}
\newcommand{\SL}{\mathrm{SL}}
\newcommand{\Q}{\mathbb{Q}}\newcommand{\A}{\adele}
\newcommand{\Zp}{\mathbb{\Z}_p}
\newcommand{\bs}{\backslash}
\newcommand{\ra}{\to}
\newcommand{\integer}{\mathrm{int}}
\newcommand{\Spec}{\mathrm{Spec}}

\newcommand{\Isom}{\mathrm{Isom}}
\newcommand{\SO}{\mathrm{SO}}
\newcommand{\X}{\mathbf{X}}
\newcommand{\K}{\tilde{K}}

\newcommand{\Galois}{\mathrm{Gal}(\overline{\Q}/\Q)}

\renewcommand{\X}{\mathbf{X}}
\newcommand{\class}{\mathscr{C}}

\renewcommand{\H}{\mathbf{H}}
\newcommand{\proj}{\Pi}

\newcommand{\Grass}{\mathrm{Grass}}

\newcommand{\codim}{\mathrm{codim}}

\newcommand{\val}{\mathrm{val}}

\newcommand{\Cliff}{\mathrm{Cliff}}
\newcommand{\Ad}{\mathrm{Ad}}
\newcommand{\G}{\mathbf{G}}

\newcommand{\etadot}{\dot{\eta}}

\newcommand{\Lie}{\mathrm{Lie}}
\newcommand{\Hom}{\mathrm{Hom}}
\newcommand{\tensor}{\otimes}

\newcommand{\xidot}{\dot{\xi}}
\renewcommand{\dim}{\mathrm{dim}}

\newcommand{\meas}{\mathrm{meas}}

\newcommand{\PGL}{\mathrm{PGL}}
\newcommand{\genus}{\mathscr{G}}
\newcommand{\spingenus}{\mathscr{G}_{\mathrm{Spin}}}
\newcommand{\adele}{\mathbb{A}}
\newcommand{\lattices}{\mathscr{L}}
\newcommand{\spinlattices}{\mathscr{L}_{\mathrm{Spin}}}
\newcommand{\order}{\mathscr{O}}

\newcommand{\ortho}{\mathrm{O}}
\newcommand{\so}{\mathrm{SO}}
\newcommand{\spin}{\mathrm{Spin}}

\newcommand{\Rspin}{\mathscr{R}_{\mathrm{Spin}}}
\newcommand{\Z}{\mathbb{Z}}
\newcommand{\beq}{\begin{displaymath}}
\newcommand{\eeq}{\end{displaymath}}
\DeclareFontFamily{OT1}{rsfs}{}
\DeclareFontShape{OT1}{rsfs}{n}{it}{<-> rsfs10}{}
\DeclareMathAlphabet{\mathscr}{OT1}{rsfs}{n}{it}
\newtheorem{lemma}{Lemma}
\newtheorem{prop}{Proposition}
\newtheorem{thm}{Theorem}
\newtheorem{defn}{Definition}

\begin{document}

\setcounter{tocdepth}{1}

\title{Local-global principles for representations of quadratic forms.}
\author{Jordan Ellenberg and Akshay Venkatesh}
\maketitle
\begin{abstract}We prove the local-global principle holds for the problem of representations of quadratic forms by quadratic forms, in codimension $\geq 7$. The proof uses the ergodic theory of $p$-adic groups, together with a fairly general observation
on the structure of orbits of an arithmetic group acting on integral points of a variety. 
\end{abstract}
\tableofcontents

\section{Introduction}

\subsection{General comments.}

Let $(\Z^n,Q)$ and $(\Z^m,Q')$ be quadratic lattices (free finitely generated
abelian groups endowed with quadratic forms.)  We say $Q'$ is {\em
  represented} by $Q$ if $(\Z^m,Q')$ can be embedded isometrically
into $(\Z^n,Q)$.  The problem of determining whether one quadratic
form represents another goes back to the beginning of modern number
theory:  for instance, Lagrange's theorem on sums of four squares says
precisely that the quadratic form $x^2 + y^2 + z^2 + w^2$ represents
every nondegenerate quadratic form of rank $1$.  The case $m=2, n=3$ (representations of binary forms by ternary forms) was already studied by Gauss in {\em Disquisitiones.}  Another question of this type (with $n=4, m=2$) is:
are there orthogonal  vectors $\mathbf{x}_1, \mathbf{x}_2$ in the
standard Euclidean lattice $\mathbb{Z}^4$ with prescribed lengths?
Schulze-Pillot's paper~\cite{schu:qfsurvey} is an excellent survey of
both classical and modern work on this problem.

We say $Q'$ is {\em everywhere locally representable} by $Q$ if the
quadratic form $Q' \otimes \Z_p$ embeds into the quadratic form $Q
\otimes \Z_p$ for every $p$, and $Q' \otimes \mathbb{R}$ embeds in $Q \otimes
\mathbb{R}$.  A result of the form ``if $Q'$ is everywhere locally
representable, it is representable'' is referred to as a {\em
local-global} principle.  Results of this kind are part of a general
program in arithmetic geometry to understand {\em Hasse principles}
for varieties:  in this case,  the representability of $Q_2$ by $Q_1$
corresponds to the existence of an integral point on a certain affine variety
$\mathbf{X}$, and such a result amounts to the statement that
$\mathbf{X}$ has an integral point if it has a $\Z_p$-point for every
$p$.

 In the general case Siegel gave a mass formula and proved a
local-global principle when $Q$ is indefinite.  In the definite case
one has an obstacle arising from the possibile nontriviality of the
genus; in other words, there may be many forms which are isomorphic to
$Q$ over every $\Zp$ and $\mathbb{R}$, but not over $\Z$.  Here we explain how to overcome this obstacle by means of the ergodic
theory of $p$-adic groups (``Ratner's theorem,'' generalized to the
$p$-adic case by Ratner \cite{Ratner} and Margulis-Tomanov \cite{MT}) and prove a local global
principle (when the minimum integer represented by $Q'$ is
sufficiently large) when $n-m \geq 7$.  The number $7$ can likely be
reduced here; it seems likely that one can achieve $n-m\geq 3$, under
certain mild ramification conditions -- such as considering only those $Q'$ whose discriminant
is not divisible by some fixed large prime $p$ -- by means
a more refined analysis of the maximal subgroups of the orthogonal group.

Previously this type of result was known -- by very different methods
-- in the range where $n \geq 2m + 3$; this result is due to Hsia,
Kitaoka, and Kneser~\cite{hsia:2m+3}.   The present method is closely
related to work of Linnik and we discuss the connections further in
\S \ref{linnik}.  In short, we are showing that a certain variety has a {\em integral point} by using ergodic theory! This aspect is quite striking to the authors and contrasts with the use of ergodic theory or dynamics to 
produce solution to Diophantine inequalities (for example, Margulis' proof of the Oppenheim conjecture).  \footnote{
Of course, it is not precisely true that ergodic theory ``produces'' an integral point; conceptually, the point
is not dissimilar to the Hardy-Littlewood method, where one deduces there exists an integral point
on an affine variety by proving equidistribution of integral points in
some larger space (an affine space in the Hardy-Littlewood setting, a
union of varieties parametrized by the genus of $Q$ in our case.). }

\subsection{Statement of theorem.} \label{sec:thmstatement}

The methods are quite robust and applicable over an arbitrary number field, and indeed our main result (Prop. 
\ref{prop:spin}) is stated in that generality, but we state the main implication only in the most classical setting.

\begin{thm} \label{thm:final}
Let $Q$ be a positive definite quadratic form on $\Z^n$. Then there exists $c:= c(Q)$ such that $Q$ represents all quadratic forms $Q'$ 
in $m \leq n-7$ variables that are everywhere locally representable, have squarefree discriminant,
and minimum $\geq c(Q)$.  
\label{th:main}
\end{thm}

We recall that the discriminant of the quadratic form $Q$ is the determinant of the matrix
$(Q(e_i + e_j) - Q(e_i) - Q(e_j))_{ij}$, where $e_i$ is a basis of $\Z^n$, and the minimum 
of $Q$ is the smallest nonzero element of $Q(\Z^n)$. 

The assertion about ``squarefree discriminant'' is probably stronger than
necessary. We note that it is possible for the local-global principle
to fail without such an assumption, as was brought to our attention by
W.K. Chan; however, we expect that one could formulate a more precise
theorem that excluded precisely such cases by a more detailed local
analysis.  Schulze-Pillot has indicated to us such a sharpened version,
utilising in particular the auxiliary condition of ``bounded imprimitivity of local representations''
(see \cite[p4]{schu:qfsurvey}, especially the paragraph after (1.8)).

 The assertion about minimum $\geq c(Q)$ (which is clearly
necessary) should be seen as a condition on local representability at
$\infty$. A defect of the method is that it does not yield an
effective upper bound for $c(Q)$.

On the other hand, the method of proof should yield a quantitative
result. This requires additional technical work and we have not aimed
for it; however, the shape of the quantitative result would be as
follows: If $r(Q,Q'),\tilde{r}(Q,Q'),g(Q)$ denote, respectively, the
weighted number of representations of $Q'$ by $Q$, the weighted number
of representations of $Q'$ by the spin genus of $Q$, and the mass of
the spin genus of $Q$, then $$r(Q,Q') \sim
\frac{\tilde{r}(Q,Q')}{g(Q)},$$ as the minimum of $Q'$ approaches
$\infty$.  Note that $\tilde{r}(Q,Q')$ can be given {\em exactly} by
the Siegel mass formula for the spin genus.

The proof of Theorem~\ref{th:main} is given in Section \ref{proof}. It is independent
of the rest of the introduction and of Section \ref{algebraicstructures}, and the 
reader interested only in this proof may proceed immediately to Section \ref{proof}. However, these
intervening sections provide (we hope) some context for the algebraic ideas underlying the method. 

\subsection{The role of the stabilizer.}\label{stabilizer}
We take a moment to describe, in a quite general context, a key feature of ``integral orbit problems'' -- i.e., problems pertaining to the orbits of an arithmetic group on the integral points of a variety -- 
utilized in this paper.   This feature has been noted by many people in many contexts in number theory.  We attempt to present a quite general (though vague) version here, and give a more precise discussion in Section \ref{algebraicstructures}.  We also refer to Section \ref{connection} for more discussion of the provenance of this type of idea. 

Let $\mathbf{G}$ be a semisimple group over $\mathbb{Q}$ that acts on
a variety $\mathbf{X}$ defined over $\mathbb{Z}$; let $\Gamma$ be a
lattice in $\mathbf{G}(\Q)$ that preserves $\mathbf{X}(\Z)$.  Then
evidently $\Gamma$ acts on $\mathbf{X}(\Z)$.  An important observation for the present paper is that the {\em the set of
orbits} $\mathbf{X}(\Z) /\Gamma$ can be described in terms of the {\em
  stabilizer} $\G_{x_0}$ of a point $x_0 \in \mathbf{X}(\Z)$.

More precisely, $\mathbf{X}(\Z)/\Gamma$ is ``closely related'' to both of the following two objects:
\begin{enumerate}
\item A fiber of the map $\Omega \backslash \G_{x_0}(\adele_f) / \G_{x_0}(\Q)
\rightarrow \Omega' \backslash \G(\adele_f)/ \G(\Q)$, where $\adele_f$
is the ring of finite adeles of $\Q$ and $\Omega,\Omega'$ are suitable open compact subgroups. 
\item A fiber of the map\footnote{Here we have used algebro-geometric language. However, this will not be used in the proof of Theorem \ref{thm:final}.} $H^1(\Spec \, \Z, \G_{x_0}) \rightarrow
  H^1(\Spec \,\Z, \G)$, where we have chosen flat models for $\G$ and
  $\G_{x_0}$ over $\Spec \, \Z$ (if this is possible) and the
  cohomology is fppf.    
\end{enumerate} 

We refer to Sec. \ref{algebraicstructures} for a full explanation of what ``closely related'' means.  
 For now we remark that the second assertion
is simply an integral version of the fact that, if $\H \subset \G$ are algebraic groups over a field $k$,
then the $\G(k)$ orbits on $\G/\H(k)$ are parameterized by the kernel of 
$H^1(k, \H) \rightarrow H^1(k,\G)$. The first remark is then not surprising, for it is well-known that
$H^1$ of algebraic groups over $\Z$ can be interpreted in terms of suitable adelic quotients.

In the context of Theorem \ref{thm:final}, we will take for $\X$ the variety parameterizing isometric embeddings
of a quadratic form $Q'$ into another quadratic form $Q$. Now, all we wish to prove
is that $\X(\Z)$ is nonempty if certain local conditions are satisfied. This will follow from establishing the {\em surjectivity}
of the maps described above. 
In the first picture, this surjectivity can be approached by studying the dynamics\footnote{At a vague level one can see the existence of this hidden dynamical structure as a refinement of the obvious fact that is possible to take $g \in \G(\Q)$ which moves one
$\Gamma$-orbit to another.} of the action
of $\G_{x_0}(\adele_f)$ on $\G_{x_0}(\Q) \backslash \G(\adele_f)$.  
In practice, there is no loss in passing from the adeles to a single completion $\Q_p$ and applying
dynamical results for actions of the $p$-adic Lie group $\G_{x_0}(\Q_p)$.

\subsection{Idea of the proof.} \label{sec:proofidea}

We now give a more concrete outline of the plan of the proof of
Theorem~\ref{th:main}.  We also highlight
historical uses of related techniques (especially by Linnik) in special cases. 
The key new ingredient (when compared to existing methods for the analysis of such problems
) is Ratner's theorem; to apply this result,
one needs in addition to reduce the question to classification of ergodic measures
-- here Lemma \ref{ergodiclimit} allows a considerable simplification
-- and to verify appropriate ``non-focussing'' conditions.


The treatment of certain auxiliary issues necessary to apply Ratner's theorem 
is deferred to the Appendix. For the purpose of the present section, this can be assumed. 
\subsubsection{Outline in elementary terms.} \label{elementary}


%

We begin with a bit of hand-waving to give the general idea.  Let
$(\Z^m,Q')$ be some quadratic form which is everywhere locally
represented by $Q$; we seek 
to prove that it is globally represented. By Hasse-Minkowski, $Q'$ is globally represented by $Q$
over the rational numbers: that is to say, there exists an isometric embedding $l_{\Q}$ of $(\Z^m, Q')$
into $(\Q^n, Q \otimes \Q)$. However, we still have a large symmetry group to play with: clearly, we can compose $l_{\Q}$ with any isometry
$\gamma \in \so_Q(\Q)$ and the result will still be an isometric embedding; thus,
we can attempt to find $\gamma$ such that $\gamma \circ l_{\Q}$ actually has image
in $\Z^n$. The idea of this paper is to use Ratner's ergodic theorems to show that one can find such a $\gamma$, and indeed with $\Q$ replaced by the much smaller ring $\Z[1/p]$ for suitable $p$. 
That the existence of such a $\gamma$ should be a rather subtle matter can already be seen in the case $n=3, m=1$, where the local-global theorem was established by Duke and Schulze-Pillot and is very closely tied to subconvexity bounds for $L$-functions. 
In any case, this description does not really capture the underlying symmetry of the situation; 
we give a more detailed description in what follows.

\subsubsection{Ternary quadratic forms and the work of Linnik.} \label{linnik}
 First, we start with a situation to which
our theorem is {\em not} applicable, but which nevertheless illustrates the main concepts that enter:
namely $n=3, m=1$, i.e. the question of representability of integers by ternary quadratic forms.
This case was essentially completely settled by W. Duke and R. Schulze-Pillot \cite{Du}; for now, we shall describe earlier
work due to Linnik that gave a weaker result \cite{linnik} but is closer to our needs. 

Gauss already observed in {\em Disquisitiones} that the number of
primitive representations $d = x^2+ y^2+ z^2$ is $12 h(-4 d)$ for $d$ congruent to $1$ mod $4$, where $h(-4d)$ is the class number of the quadratic order
$\mathscr{O}_d := \Z[\sqrt{-d}]$.  From now on we assume $d$
squarefree to avoid having to repeatedly specify that we consider only
primitive representations.  Gauss's formula can easily be
understood in the framework of Section \ref{stabilizer}; here the
group $\G = \mathrm{SO}(3)$, the space $\X$ is the quadric
$x^2+y^2+z^2 = d$, and the stabilizer $\G_{x_0}$ is a form of 
$\mathrm{SO}(2)$.

One way to interpret Gauss' formula is to 
construct an explicit map from solutions
$\{(x,y,z) \in \mathbb{Z}^3: x^2+y^2 + z^2 =d \} $
to binary quadratic forms of discriminant $-4d$; such a map is given by associating
to a solution $(x,y,z)$ the restriction of the Euclidean quadratic form to the orthogonal
complement $(x,y,z)^{\perp}$. 
Although this is not a bijection, one can precisly quantify how far it is from being a bijection. 

However, there is a more suggestive (for our purposes) way of phrasing the answer, that is more familiar
from the context of Heegner points on division algebras. The set 
\begin{equation} \label{solutions} \{(x,y,z) \in \mathbb{Z}^3: x^2+y^2+z^2 = d\}/\mathrm{SO}_3(\Z), \end{equation}
  carries an action of $\mathrm{Pic}(\order_{d})$. This action is {\em almost} simply transitive; it fails to be simply transitive because of problems at $2$. More precisely, it is transitive, and its kernel
 is the $2$-torsion ideal class generated by a prime ideal of $\Z[\sqrt{-d}]$ above $2$. 
This Picard group is a quotient of the idele class group of $\Q(\sqrt{-d})$; in particular, if $p$
is a split prime in $\Q(\sqrt{-d})$, then the group $\mathbb{Q}_p^{\times}$ acts
on (\ref{solutions}); this action is trivial on $\Z_p^{\times}$
and indeed factors through a finite cyclic group.  Such actions are discussed
in a more general context in Section \ref{algebraicstructures}. 


This particular case is not relevant to our paper, since the genus of the quadratic form
$x^2+y^2+z^2$ contains only one element, and so the local-global principle is evident. 
If one had replaced $x^2+y^2+z^2$ by a general definite ternary quadratic form $Q(x,y,z)$, then 
the relevant algebraic statement is the following: Let $\{Q_1, \dots, Q_g\}$ be the genus of $Q$. 
Then
$$\bigcup_{1 \leq i\leq g} \{(x,y,z) \in \mathbb{Z}^3: Q_i(x,y,z) = d\}/\mathrm{SO}_{Q_i}(\Z)$$
still is a principal homogeneous space (or ``almost'' a principal homogeneous space) for a suitable Picard group and carries
an action of $\Q_p^{\times}$ for any $p$ that is split in $\Q(\sqrt{-d D})$, where $D$ is the discriminant of $Q$.  
Linnik's method
(in modern language) is then to interpret this action in terms of a suitable collection of closed orbits of $\Q_p^{\times}$
on a $p$-adic homogeneous space (i.e., the quotient of a $p$-adic Lie group by a lattice),
and then to prove equidistribution results about this collection of closed orbits.

A modern interpretation and extension of Linnik's work will appear in the second paper of the sequence \cite{ELMV}, 
and further work along these lines will appear in \cite{MV}.

\subsubsection{Higher rank quadratic forms and class number problems.}
 A key observation of this paper is that, in the higher rank case, one retains a residue of this type of 
structure after passing to a suitable covering set; however, rather than the action of the $p$-adic torus $\Q_p^{\times}$ (which one can think of as $\mathrm{SO}_2(\Q_p)$), one obtains
the action of a special orthogonal $p$-adic group in more variables.  Again, this can be understood
in terms of Section \ref{stabilizer}; we first describe the action in more classical terms and then give a ``dictionary'' between this description and Section \ref{stabilizer}.

Let $Q$ be a positive definite quadratic form of rank $n$, $Q'$ a positive definite form of rank $m$.

Let $\mathscr{G} = \{Q=Q_1, Q_2, \dots, Q_g\}$ be the genus of $Q$.
Let $\mathscr{R}$ be the set of isometric embeddings of the lattice
$(\Z^m, Q')$ into any of the lattices $(\Z^n, Q_i)$ for some $1 \leq i
\leq g$.  We should like to know whether the forgetful map $\mathscr{R} \to
\mathscr{G}$ is surjective; Siegel's mass formula gives explicit
formulas for the size of $\mathscr{R}$ and $\mathscr{G}$
(appropriately weighted), but {\em gives no information about the
nature of the map from $\mathscr{R}$ to $\mathscr{G}$.}

Unlike the case $n=3, m=1$ there are no {\em direct} group actions on
$\mathscr{R}$; the action in that case was a special feature arising
from the fact that the stabilizer (a form of $\mathrm{SO}_2$) was
abelian.  What remains true in general is that we can cover $\mathscr{R},
\mathscr{G}$ by {\em profinite} sets $\tilde{\mathscr{R}},
\tilde{\mathscr{G}}$:

\begin{equation} \begin{CD}
\tilde{\mathscr{R}} @>>> \tilde{\mathscr{G}} \\  @VVV @VVV \\
\mathscr{R} @>>> \mathscr{G} 
  \end{CD} \end{equation}

Here:
\begin{itemize}
\item 
$\tilde{\mathscr{G}} = \Gamma \backslash G$, where $G$ is a $p$-adic spin group in $n$ variables and $\Gamma$ is a lattice in $G$;
\item The image of $\tilde{\mathscr{R}}$ in $\tilde{\mathscr{G}}$ is an orbit of a certain subgroup
$H \subset G$, a spin group in $n-m$ variables.
\end{itemize}

The $p$-adic version of Ratner's theorem allows us to understand that
this $H$-orbit is (fairly) dense in $\tilde{\mathscr{G}}$, so
$\mathscr{R} \to \mathscr{G}$ is surjective.  More precisely, we show
that every open subset of $\tilde{\mathscr{G}}$ (in particular, the
preimage of an element of $\mathscr{G}$) has nontrivial intersection
with all but finitely many of the $H$-orbits that arise (for various $Q'$) in the above discussion. 

The dictionary between the discussion above and Section
\ref{stabilizer} is as follows: we take $\G = \so(Q)$ and $\X$ to be
the variety parameterizing isometric embeddings of $Q'$ into $Q$, i.e.
the variety of linear maps $\ell: \Z^m \rightarrow \Z^n$ so that $Q
\circ \ell = Q'$.  The stabilizer $\G_{x_0}$ of a point $x_0 \in
\X(\Q)$ is then an orthogonal group in $n-m$ variables. Then
$\mathscr{R}$ (resp. $\tilde{\mathscr{R}}$) corresponds to $\Omega
\backslash \G_{x_0}(\adele_f) / \G_{x_0}(\Q)$
(resp. $\G_{x_0}(\adele_f)/\G_{x_0}(\Q)$) whereas $\mathscr{G}$
(resp. $\tilde{\mathscr{G}}$) corresponds to $ \Omega' \backslash
\G(\adele_f)/ \G(\Q)$ (resp. $\G(\adele_f)/\G(\Q)$).  As we remarked
before, it is possible to replace the role of $\adele_f$ by $\Q_p$ for
a suitable $p$.

\subsection{Connection to existing work.} \label{connection}
Schulze-Pillot has pointed out to us that the set-up of the proof of Theorem \ref{thm:final} is quite close
to that of Hsia, Kitaoka and Kneser \cite{hsia:2m+3}. In essence, when the proof of Theorem \ref{thm:final} is unwound, we pass to the ring $\Z[1/p]$, i.e.
allow denominators at a suitable auxiliary prime $p$, and then pass back to $\Z$ (cf. description in Section \ref{elementary}). This is also done in \cite{hsia:2m+3}. 

As far as the ergodic side of the present paper goes, the closest cognate to our work is in the paper \cite{EO} of Eskin and Oh. They consider
a situation analogous to that discussed in the first situation of Sec. \ref{stabilizer}
but when the stabilizer $\mathbf{G}_x$ has noncompact real points. In that case, there is no issue of local-global principle (for, in the cases considered, the stabilizer $\mathbf{G}_x$ is semsimple
and satisfies a suitable version of strong approximation); the concern of \cite{EO} is instead
to prove uniform distribution results for integral points, using, in that case, the results of Ratner for real groups
and the results of Dani-Margulis \cite{DM}. In our ($p$-adic) setting, we do not have the results
of \cite{DM} available; Appendix A gives a self-contained proof (assuming the classification
of ergodic measures \cite{MT}, \cite{Ratner}) of what we need.

As for the arithmetical side of the present paper, the presence of the kind of group actions remarked on in Section \ref{stabilizer} and elaborated in Section \ref{algebraicstructures}
has been noted in many different instances, though perhaps not in a unified way. We mention in particular 
the work of Linnik \cite{linnik}, Eichler \cite{eichler}, Kneser \cite{Kneser} and Weil \cite{weil} on quadratic forms; the latter two papers already contain a framework essentially equivalent to what we use for the proof of Theorem \ref{thm:final}. More recently, we refer to the 
work of Shimura \cite{shimura} and the beautiful results of Bhargava \cite{Bhargava}, which
of course go much deeper.

Let us briefly contrast the present work with results on Diophantine inequalities. 
When considering an irrational quadratic forms in $n$ variables from the point of view of Diophantine {\em inequalities}, e.g. the Oppenheim conjecture, it is natural to consider the action and dynamics
of an orthogonal group in $n$ variables. However, when investigating the arithmetic properties
of a rational quadratic form in $n$ variables, we will be led naturally to consider the action
of a $p$-adic $O(n-1)$. In general, a fundamental difference between irrational and rational quadratic forms seems to be the following: whereas for irrational quadratic forms one may utilize
only the dynamics of a {\em real} orthogonal group, one may study rational quadratic forms 
through the dynamics of an {\em adelic} orthogonal group. This added freedom is precisely why 
we are able to say something about positive definite quadratic forms: though the real points
of the associated orthogonal group are compact, the $p$-adic points need not be.

We must emphasize that from the ergodic point of view there is not much novelty except, perhaps, 
Lemma \ref{ergodiclimit} in the appendix. The ``deep'' and important ingredient
is the classification of ergodic measures, due to Ratner and Margulis/Tomanov in the setting we consider. 

\subsection{Acknowledgements.}
We would like to express, first of all, our gratitude to W.K. Chan,  Hee Oh and Rainer Schulze-Pillot: all of whom took the time and trouble of reading parts of the present paper quite carefully and
sent us very helpful comments and mathematical corrections.   The exposition (and also correctness!)
of the sections concerning ergodic theory has been considerably improved through Oh's comments;
similarly, the exposition of the sections concerning the arithmetic of quadratic forms have been improved
greatly through the comments of Chan and Schulze-Pillot.  Schulze-Pillot also provided many references  to put the present work in the correct context relative to the quadratic forms literature. 

We would also like to thank Manfred Einsiedler and Elon Lindenstrauss for many discussions about ergodic theory, Jon Hanke and Peter Sarnak for their help with quadratic forms, and James Parson
for  insightful letters concerning Section \ref{algebraicstructures}. 

The first author was partially supported by NSF-CAREER Grant DMS-0448750 and a Sloan Research
    Fellowship; the second author was supported by a Clay Math Research Fellowship and NSF Grant
    DMS-0245606. We thank the Clay Mathematics Institute for supporting collaborative visits during which this paper was written.

\section{Proof of Theorem \ref{thm:final}} \label{proof}

The scheme of proof is as follows.  In Section \ref{QS} we give some background on quadratic spaces
over global and local fields. In Section \ref{SGR} we introduce the notion
of {\em spin globally representable} and state the main Proposition \ref{prop:spin}
which is valid over an arbitrary number field. In Section \ref{implies}, we show that Proposition \ref{prop:spin}
implies Theorem \ref{thm:final}. In Section \ref{sec:reduction} and Section \ref{Ratner}, we explain how Proposition \ref{prop:spin} is reduced to a statement which can be approached by Ratner's theorem, together with a
result about generation of spin groups by embedded spin groups of smaller dimension. Finally, in Section \ref{sec:grouptheory}
we resolve the necessary group-theoretic issues, concluding the proof. 

\subsection{Quadratic spaces, lattices, genera.} \label{QS}
We begin with some relatively standard material on quadratic spaces. 

Let $F$ be a number field, $\order$ the ring of integers of $F$. 
Let $(V,q)$ be a quadratic space over  $F$.  By a {\em
  lattice} in $V$ we mean a locally free $\order$-submodule of $V$
whose rank is $\dim V$.  Let $\Lambda_V \subset V$ be a lattice on which $q$ is integral, i.e., a lattice such that $q(\Lambda_V) \subset \order$.

We shall assume that $q$ is definite at all infinite places of $F$. \footnote{Otherwise, the spinor genus
of $\Lambda_V$ contains a unique class by strong approximation as soon
as $\dim(V) \geq 3$.  In that case we shall regard the problem as solved, although in the general (nonsquarefree discriminant case) there are complicated local issues involved. Note that the assumption entails that $F$ is totally real, but we shall not make any use of this.}
Attached to $q$ we have a bilinear form $\langle v_1, v_2 \rangle := \frac{1}{2} (q(v_1+v_2) - q(v_1) - q(v_2))$. This bilinear form is not necessarily integral, but takes values in $\frac{1}{2} \order$.

Let $\GL(V) , \ortho_V,\so_V,\spin_V$ be (respectively) the general linear, orthogonal, special orthogonal,
and spin groups of $V$.  These are algebraic groups over $F$ and consequently we may speak
of their points over any ring containing $F$. If $\delta \in V$ is such that $q(\delta) \neq 0$
we will denote by $r_{\delta} \in \ortho_V(F)$ the reflection through the orthogonal complement of $\delta$: that is to say $w \mapsto w - 2 \frac{\langle w, \delta \rangle}{\langle \delta, \delta \rangle} \delta$.

Let $\adele$ (resp. $\adele_f$)
be the ring of adeles (resp. finite adeles) of $F$.

 It is well-known that $\GL(V, \adele_f)$ acts on the
lattices in $V$; by restriction we obtain an action of $\ortho_V(\adele_f)$ on 
the lattices in $V$; via the map $\spin_V(\adele_f) \to \ortho_V(\adele_f)$, 
we obtain also an action of $\spin_V(\adele_f)$ on lattices. 
Recall that one says that two lattices $\Delta_1, \Delta_2$ in $V$ are {\em locally isomorphic}
if they are isomorphic as quadratic spaces over each completion of $\mathscr{O}$. 
With this definition, the equivalence relation corresponding to the $\ortho_V(\adele_f)$-orbits
is exactly that of local isomorphism. 

If $L_1, L_2$ are two locally free $\order$-modules endowed with quadratic forms, we denote by
 $\Isom(L_1, L_2)$ the set of isometric embeddings of $L_1$ into $L_2$.

For each finite place $v$, the stabilizer of a lattice $\Delta$ in 
$\so_V(F_v)$ is an open compact subgroup $K_{\Delta,v}$. 
Let $\K_{\Delta,v}$ be the preimage of $K_{\Delta,v}$ in 
$\spin_V(F_v)$. 
Put $K_{\Delta,f} = \prod_{v \, \mathrm{finite}} K_{\Delta,v}$
and $\K_{\Delta,f} = \prod_{v \, \mathrm{finite}} \K_{\Delta,v}$. 
 In the case $\Delta = \Lambda_V$ we write simply $K_v, K_f, \K_v, \K_f$.
%
%

We recall that, for each place $v$, one has a homomorphism (the ``spinor norm'')
from $\so_V(F_v)$ to $F_v^*/(F_v^*)^2$, which sends the product of reflections $r_v r_{v'}$
to $q(v) q(v')$. Moreover the image of $\spin_V(F_v)$ in $\so_V(F_v)$ coincides
with the kernel of the spinor norm.

We will need some facts about quadratic forms over local fields. Continue to assume that
$v$ is a finite place; let $\order_v$ be the closure of $\order$ in $F_v$.

\begin{lemma} \label{solvability}
Suppose $J$ is a nondegenerate quadratic space over $F_v$, and the residue characteristic
of $F_v$ is larger than $2$.
Then:
\begin{enumerate}

\item If $\dim(J) \geq 3$, then the
spinor norm $\so_J(F_v) \to F_v^*/(F_v^*)^2$ is surjective.

\item If $\dim(J) \geq 5$, then $J$ is isotropic. 
\item If $\dim(J) \geq 5$, the spin group $\spin_J(F_v)$ is generated by the 
unipotent radicals of parabolic subgroups. 
\item If $\dim(J) \geq 5$, then the spin group $\spin_J(F_v)$ is generated by the embedded
spin groups $\spin_P(F_v)$ for $P$ a hyperbolic plane\footnote{Recall that a hyperbolic plane
is a two-dimensional nondegenerate quadratic space possessing an isotropic vector.}  inside $J$. 
\item If $\dim(J) \geq 5$, then the subgroup of $O_J(F_v)$ generated by reflections associated to vectors of length $1$
contains the image of $\spin_J(F_v)$ in $O_J(F_v)$. 
\end{enumerate}
\end{lemma}
\proof
We verify first the assertion about the surjectivity of the spinor norm.
We may assume $\dim(J)=3$.  The question is unchanged by replacing the quadratic
form on $J$ by any multiple of itself; thus, in suitable coordinates, the form takes the shape
$x^2 +q(y,z)$, where $q(y,z)$ is a quadratic form in $y,z$. If $q$ represents the value $d \in F_v^*$,
then it easily follows that all norms from the quadratic extension $F_v(\sqrt{-d})$ are values
of the spinor norm on $\so_J(F_v)$. By class field theory, if $d$ is not a square,  the group of norms
from $F_v(\sqrt{-d})$ is an index $2$ subgroup of $F_v^*$ determining the square class of $d$. 
So it suffices to show that the nonzero values taken by $q$ do not all lie within a single coset
of $F_v^*/(F_v^*)^2$. But $q$ itself is a multiple of a norm form on a quadratic extension 
of $F_v$, whence the assertion. 
%


%

We omit the proof of the second assertion, which is due to Hasse.

If $J$ is isotropic, the group $\spin_J(F_v)$ is projectively simple
(any normal subgroup is central, and in particular finite), as is proved in 
\cite[Theorem 5.27]{Artin}. This implies the third, fourth and fifth assertions.
%
\qed

Let $v :  F_v \to \mathbb{Z} \cup \infty$ be the associated valuation. 
Recall that we say that a subspace of a quadratic space is {\em regular} if the induced quadratic
form is nondegenerate. 
For any quadratic
subspace $Z$ of $V \otimes F_v$ we put
$\val(Z)$ to be the valuation of the discriminant of the quadratic form induced
on $Z \cap (\Lambda \otimes \order_v)$. In other words, choosing an $\order_v$
basis $w_1, \dots, w_r$ for $Z \cap (\Lambda \otimes \order_v)$, 
we put $\val(Z) := v(\det(2 \langle w_i, w_j \rangle))$.  (The inclusion of $2$ is to guarantee
that $\val(Z) \geq 0$, and is superfluous
if the residue of characteristic of $F_v$ is bigger than $2$, as will always the case in our discussion.)
$Z$ is regular if and only if $\val(Z) < \infty$.



%

\begin{lemma} \label{lem:compact}
Suppose $Z_i$ is a sequence of subspaces of $V \otimes F_v$
such that $\val(Z_i)$ remains bounded. 
Then there exists a compact set $\Omega \subset \spin_V(F_v)$
and a partition of $(Z_i)$ into finitely many subsequences,
such that if $Z_i, Z_j$ belong to the same subsequence there exists $\omega_{ij} \in \Omega$
with $\omega_{ij} Z_i = Z_j$. 
\end{lemma} \proof
Without loss, we may assume that $\dim(Z_i)$ is fixed, say $=r$. Let $\Grass_r$ 
be the Grassmannian of $r$-dimensional subspaces in $V$. Then
the result follows easily from the following assertions:
\begin{enumerate}
\item For any point $Z \in \Grass_r(F_v)$, the map
$\mathrm{GL}_V \mapsto \Grass_r$ given by $g \mapsto g. Z$ is submersive\footnote{We say a map of two smooth algebraic varieties $\mathbf{V}_1 \to \mathbf{V}_2$ over a field $k$ is submersive at $v_1 \in \mathbf{V}_1(k)$
 if the induced map on tangent spaces is surjective. If $k$ is a local field this implies (``implicit function theorem'') that the image
 contains a neighbourhood of the image of $v_1$.}
at the identity (in particular, the image of an open neighbourhood of the identity in $\GL_V(F_v)$ contains an open neighbourhood of $Z$). 
\item For any point $Z \in \Grass_r$ parametrizing a {\em regular} subspace, the map
$\spin_V(F_v) \mapsto \Grass_{r}$ given by $g \mapsto g. Z$ is submersive  
at the identity (in particular, its image of an open neighbourhood of the identity contains an open neighbourhood of $Z$). 
\item The set $\{Z \in \Grass_{r}: \val(Z) \leq N\}$ is compact. 
\end{enumerate}
The first two assertions may be checked at the level of tangent spaces.  
For the final assertion
it suffices to check that the complement of the subset in question is {\em open}. 
This follows from the fact that $g \mapsto q_v^{-\val(g.Z)}$ is {\em continuous} 
for $g \in \GL_V(F_v)$. \footnote{Indeed, choose any $g_0 \in \GL_V(F_v)$; we claim
there is a neighbourhood on which $g \mapsto \val(g \cdot Z)$ is the valuation
of a polynomial function in the coordinates of $g$. Let $U_v \subset \GL_V(F_v)$
be the stabilizer of $\Lambda \otimes \order_v$. For $u \in U$, we have $u g_0 Z \cap (\Lambda \otimes \order_v) =
u (g_0 Z \cap (\Lambda \otimes \order_v))$; thus, the map $u \mapsto \val(u g_0 Z)$
is the valuation of a polynomial in the coordinates of $u$. }
 \qed

Let $\lattices(V) = \so_V(\adele_f) \cdot \Lambda_V$ be the set of lattices in $V$
that are locally isomorphic to $\Lambda_V$,
$\spinlattices(V) = \spin_V(\adele_f) \cdot \Lambda_V$ the set of lattices
``locally spin-isomorphic'' to $\Lambda_V$.
Then $\lattices(V)$ is identified
with the quotient $\so_V(\adele_f)/K_f$, and $\spinlattices(V)$ is identified with $ \spin_V(\adele_f)/\K_f$. 
The {\em genus}  $\genus(V)$ of $\Lambda_V$ is the quotient of $\lattices(V)$ 
 by $\so_V(F)$. The {\em spin genus} $\spingenus(V)$ of $\Lambda_V$ is
 the image of $\spinlattices(V)$ in $\genus(V)$. 
It is well-known that 
$\genus(V)$ and so also $\spingenus(V)$ are finite sets. 

Moreover, $\genus(V)$ is identified with $\so_V(F) \backslash \so_V(\adele_f)/K_f$; moreover, 
 if we write, for each $v$, $\Theta_v$ for the image of $\spin_V(F_v) \to \so_V(F_v)$,
 and put $\Theta_f = \prod_{v \, \mathrm{finite}} \Theta_v$,  
 then $\spingenus(V)$ is identified with $\so_V(F) \backslash \so_V(F) \Theta_f K_f/K_f$. 

\subsection{The notion of spin globally representable.} \label{SGR}

Let $W$ be a
regular subspace of $V$ over $F$, with induced quadratic form $q_W$, and $\Lambda_W  = W \cap \Lambda_V$ the induced lattice. 
Our main concern in the present document is to show that $\Lambda_W$, endowed
with the quadratic form obtained from $q_W$, embeds primitively isometrically into every lattice in the spin genus of $\Lambda_V$. (We say an embedding $\ell: \Lambda_W \to \Lambda'$
is {\em primitive} if the image of $\Lambda_W$ is saturated in $\Lambda'$, i.e.
$\ell(\Lambda_W) . F \cap \Lambda' = \ell(\Lambda_W)$.  If
$disc(\Lambda_W)$ is squarefree, which will be the case for us, then any embedding $\Lambda_W \ra
\Lambda'$ is automatically primitive.) 
If this is the case, we shall say that $W$ is {\em spin globally representable}.  

In other words: a subspace $W$ is spin globally representable if, for every
$g \in \spin_V(\adele_f)$, there exists a primitive isometric embedding of the lattice $W \cap \Lambda_V$  into
$g. \Lambda_V$. 

Fix a nonarchimedean place $w$ of $F$, with residue characteristic $>2$. 
We shall say a subspace $W \subset V$ is {\em good} if 
$\codim(W) \geq 7$ and the $w$-valuation of the determinant of $\Lambda_W$
is $\leq 1$, i.e. $\val(W \otimes F_w) \leq 1$ in the notation introduced prior to Lemma \ref{lem:compact}.
There is surely scope for considerable relaxation of both these conditions.

%
%

 The following result (which we state over a general number field) is our key result, and implies
 almost immediately Theorem \ref{thm:final}. 
 \begin{prop}\label{prop:spin}
 There exist a finite list of nontrivial subspaces $E_1, \dots, E_d \subset V$
 such that any good subspace $W$ that does not contain any $E_i$
 is spin globally representable. 
\end{prop}

\subsection{Proposition \ref{prop:spin} implies Theorem \ref{thm:final}} \label{implies}
We explicate how this Proposition implies Theorem \ref{thm:final}. 

We specialize to the field $F = \Q$ and will use the classical language of quadratic forms. Let $(\Z^n, Q)$ be a positive definite quadratic form on $\Z^n$. Let $Q_1, \dots, Q_g$ be the spin genus of $Q$.  (In the language of Section \ref{QS}, with $(V= \Q^n, q=Q)$,
the $Q_i$ are the quadratic forms induced
on a set of representatives for $\spingenus(V)$.)

Let $Q'$ be a quadratic form with squarefree discriminant on $\Z^{m}$
which is everywhere locally represented by $Q$.
In what follows, we sometimes write $Q_R$ for the quadratic form induced by $Q$
on $R^n$, for an arbitrary ring $R$. 

\begin{lemma}\label{HsiaLem}
$Q'$ is globally represented by a form in the spin genus  $\{Q_1, \dots, Q_g\}$.
\end{lemma}

See \cite{Hsia}, which proves a slightly stronger  assertion, {\em without} any assumption of squarefree discriminant on $Q'$.  We include a proof in the interest of keeping the paper self-contained. 

\proof While this may be proven with the mass formula, we prefer to give a direct proof.
By Hasse-Minkowski, we may choose\footnote{This ``subspace version'' of Hasse-Minkowski is easily decuced from the usual version; see \cite[Theorem 66:3]{OMeara}} a subspace $W \subset \Q^n$
such that the restriction of $Q_{\Q}$ to $W$ is isomorphic to $Q'_{\Q}$. 
Choose $L' \subset W$ so that  the quadratic form induced on $L'$
is isomorphic to $Q'$.
Moreover, by the definition of ``local representability'' we have for each $p$ a subgroup $L_p' \subset \Z_p^n$
such that the restriction of $Q_{\Z_p}$ to $L_p'$ is isomorphic to $Q'_{\Z_p}$.

We may choose for each $p$, an element $g_p \in \spin_{Q}(\Q_p)$ with the property that
$g_p L_p' = L' \otimes \Z_p
$, in such a fashion that $g_p$ fixes $\Z_p^n$ for almost all $p$. 
Indeed there exists an isometry $h: L_p' \to L' \otimes \Z_p$
of quadratic $\Z_p$-modules; extend $h \otimes \Q_p$ to a global isometry by Witt's theorem.
This gives an element $g_p \in \so_Q(\Q_p)$ with the property that $g_p L_p' = L' \otimes \Z_p$.  To do better, we just note 
that because $\dim(T_p) \geq 3$, where $T_p$ is the orthogonal complement of $W_p$, 
the spinor norm on $\so_{T_p}$ is surjective by Lemma
\ref{solvability}. So we can modify $g_p$ by an element of $\so_{T_p}$,
thought of as an element of $\so_Q$ stabilizing $W_p$, to be in $\spin_Q(\Q_p)$. 

We set $\Lambda_p = g_p \Z_p^n$
and $\Lambda = \{\lambda \in V: \lambda \in \Lambda_p \mbox{ for all } p\}$. 
Then $(\Lambda, Q|_{\Lambda})$ is evidently in the spin genus of $(\Z^n, Q)$.

We also note that $W \otimes \Q_p = L' \otimes \Q_p = (g_p L'_p) \otimes \Q_p$, so that $(\Lambda \cap W)_p = g_p( \Z_p^n \cap L_p' \otimes \Q_p) 
= g_p L_p'$, where at the last stage we have used the fact that the
discriminant of $Q'$ is squarefree. 
So $\Lambda \cap W = L'$, that is to say, $Q'$ is represented by the quadratic form $Q|_{\Lambda}$
which belongs to the spin genus of $(\Z^n, Q)$.
\qed

We may now complete the proof of Thm. \ref{thm:final}. 
\proof (of Thm. \ref{thm:final})  Let $Q$ be as in the statement of the theorem. Let $n -m \geq 7$ and let $Q'_i$, for $i\geq 1$,
be a sequence of quadratic forms with squarefree discriminant on $\Z^{m}$,
with minima approaching $\infty$,  and all of which are locally representable by $Q$.
By the previous Lemma (or the mass formula), $Q'_i$ is represented by a form in the spin genus of $Q$.

Let $\{Q=Q_1, \dots, Q_g\}$ be the spin genus of $Q$.
Partitioning $(Q_i')$ into subsequences, we may assume
that all the $Q_i'$ embed into a fixed $Q_j$, say $Q_h$ for some $1 \leq h \leq g$.
Realize $Q_h$ as a quadratic form on $\Z^n$. So,
 in other words, we have submodules $L_i \subset \Z^n$ such
that the quadratic form induced by $Q_h$ on $L_i$ is isomorphic to $Q'_i$.
Because the discriminant of each $Q'_i$ is squarefree, we have
$L_i = \Z^n \cap W_i$, with $W_i = \Q. L_i$.

Fix a prime $p$; we will apply Proposition \ref{prop:spin} with $w=p$, $F=\Q, V = \Q^n, \Lambda_V = \Z^n$.  
By that Proposition, there is a finite collection of nontrivial subspaces $\{E_1, \dots, E_d\} \subset \Q^n$ such that any good $W$ not containing any $E_j$ is spin globally representable. 
The $W_i = \Q. L_i$ are automatically good in the sense defined prior to Proposition \ref{prop:spin} (because of the assumption of squarefree discriminant
and of codimension $\geq 7$). Moreover,
the $W_i$ cannot contain any $E_j$ if $i$ is large enough; for otherwise the minimum of $Q_i'$
would not approach $\infty$. 

%

Applying Prop. \ref{prop:spin}, we conclude that the $W_i = \Q. L_i$ are spin globally representable 
for sufficiently large $i$. 
Translating back to quadratic forms, this means precisely that $Q_i'$ embeds into {\em each} $Q_j$, for $1 \leq j \leq g$
and sufficiently large $i$; in particular, all but finitely many $Q_i'$ are represented by $Q_1=Q$. So we are done. 
\qed

 \subsection{Reduction of Prop. \ref{prop:spin} to Ratner's theorem.} \label{sec:reduction}
Our aim is now to prove Proposition \ref{prop:spin}.  

In the setting of Prop. \ref{prop:spin}, 
let $T$ be the orthogonal complement of $W$ (since $W$ is regular,
we have $W \oplus T = V$) and define
$\GL(T), \ortho_T, \so_T, \spin_T$ accordingly.
These groups are embedded in $\GL(V), \ortho_V, \so_V, \spin_V$ respectively,
and, in this embedding, they are identified with the subgroups that fix $W$ pointwise. 
\footnote{Indeed, there is a natural map of Clifford algebras $\Cliff(T) \to \Cliff(V)$,
which is injective.  This induces an injective map $\spin_T \to \spin_V$; clearly the image
is contained in the subgroup fixing $W$ pointwise. To see that this is in fact the image,
it suffices to check that $\Cliff^{\mathrm{even}}(T)$  is exactly the centralizer of $W$ in $\Cliff^{\mathrm{even}}(V)$.
Take an orthonormal basis $\{e_1, \dots, e_r\}$ for $T$ and extend it to an orthonormal basis
$\{e_1, \dots, e_n\}$ for $V$. For a subset $J \subset \{1,\dots, r\}$, we put
$e_J = \prod_{i \in J} e_i$, where the product is taken w.r.t. an increasing ordering of the elements in $J$. Then it is easy to verify that $e_j e_J e_j^{-1} = (-1)^{|J|} e_J$ if $j \notin J$
and $(-1)^{|J|+1}$ otherwise. The assertion follows from this. }

Consider the set $\Rspin(W,V)$ of lattices $\Delta \in\spin_V(\adele_f) \cdot \Lambda_V$ with the property that 
$\Delta \cap W = \Lambda_W$. Then the action of $\spin_T(\adele_f)$ preserves $\Rspin(W,V)$.
There is a natural map $\Rspin(W,V) \mapsto \spingenus(V)$, namely, that
which sends a lattice $\Delta$ to its class $[\Delta]$ in the spin genus.
Moreover, $W$ is spin globally representable if this map is surjective. 
Note that $\Lambda_V \in \Rspin(W,V)$, by definition, and so also
$\spin_T(\adele_f) \cdot \Lambda_V \subset \Rspin(W,V)$. 
To show that $\Lambda_W$ embeds primitively isometrically into every lattice in the spin genus of $V$, it will suffice,
then, to show that  $\spin_T(\adele_f) \cdot \Lambda_V \subset \spinlattices(V)$
surjects onto $\spingenus(V)$. 
For this, it will suffice to check that the closed subset 
$\spin_T(F) \backslash \spin_T(\adele_f) $
of  $\spin_V(F) \backslash \spin_V(\adele_f)$
intersects each $\K_f$-orbit.

We will prove Prop. \ref{prop:spin} in the following formulation. As in that Proposition,
we regard as fixed a certain nonarchimedean place $w$ of $F$, with residue characteristic $> 2$;
the notion of {\em good} is defined w.r.t. this place. 

\begin{prop} \label{prop:spinbis}
Let $W_i \subset V$ be a sequence of good subspaces, with the property that no infinite subsequence
of the $W_i$ has a nontrivial common intersection. Then $W_i$ is spin globally representable for all sufficiently large $i$.
\end{prop}

In fact, it is clear that Prop. \ref{prop:spin} implies
Prop. \ref{prop:spinbis}.  We now explain how Prop. \ref{prop:spinbis}
implies Prop. \ref{prop:spin}.  Suppose that Prop.~\ref{prop:spin} is
false.  We define a sequence $W_1, W_2, \ldots$ of subspaces of $V$
inductively as follows: let $\Sigma_k$ be the set of nonempty
intersections of subsets of $W_1, \ldots, W_k$, and let $W_k$ be a
subspace of $V$ which is not spin globally representable and does not
contain any subspace in $\Sigma_{k-1}$.  (We can chose such a subspace
by the negation of Prop~\ref{prop:spin}.)  It is then clear that any $r$-fold intersection among the $W_i$ has codimension at least $r$.  In particular, no infinite subsequence of $W_1,
\ldots, $ has nontrivial intersection, contradicting
Prop. \ref{prop:spinbis}.




In the setting of Prop. \ref{prop:spinbis}, we have now a sequence of orthogonal complements $T_i$ and
 associated groups $\ortho_{T_i}, \so_{T_i}, \spin_{T_i}$ etc.

Prior to applying Ratner's theorem, we must switch from an adelic to an $S$-arithmetic setting. We recall we have fixed a nonarchimedean place $w$ of $F$.
Set $\K^{(w)} = \prod_{v \neq w} \K_v$, where the product is restricted to finite places. 

By the strong approximation theorem, 
it follows that $\spin_V(F) \cdot \spin_V(F_w) \cdot \K_f = \spin_V(\adele_f)$. 
Then the adelic quotient $\spin_V(F) \backslash \spin_V(\adele_f)/ \K^{(w)}$ is naturally identified with 
$\Gamma \backslash \spin_V(F_w)$, where $\Gamma$ is the projection of
$\spin_V(F) \cap \spin_V(F_w) \tilde{K}^{(w)}$ to $\spin_V(F_w)$. Note that $\Gamma$
is a cocompact lattice in $\spin_V(F_w)$. 

We want to show that -- at least for big enough $i$ -- the quotient $\spin_{T_i}(F) \backslash \spin_{T_i}(\adele_f)$ intersects
each $\K_f$-coset in $\spin_V(F) \backslash \spin_V(\adele_f)$; it will suffice
to see that $\Gamma \cap \spin_{T_i}(F_w) \backslash \spin_{T_i}(F_w)$ intersects
each $\K_w$-coset in $\Gamma \backslash \spin_V(F_w)$.

We now wish to reduce to a situation where we are studying orbits of a {\em fixed} group (not a
varying sequence of groups like the $\spin_{T_i}(F_w)$).

Note that $W_i  \otimes F_w$ is a certain quadratic subspace of $V \otimes F_w$;
in view of the assumption that $W_i$ is good, 
$\val(W_i \otimes F_w) \leq 1$ and Lem. \ref{lem:compact} is applicable. 
Partitioning our original sequence $(W_i)$ into appropriate
subsequences, 
we may assume that there is a fixed $W_w \subset V \otimes F_w$ and a compact subset $\Omega \subset \spin_V(F_w)$ such that,
for each $i$, there is $\xi_i \in \Omega$ with
$\xi_i(W_w) = W_i \otimes F_w$.

Let $T_w$ be the orthogonal complement of $W_w$ inside $V \otimes F_w$,
and denote by $\spin_{T_w}$ the spin group of $T_w$.
Then $\spin_{T_i}(F_w) = \xi_i \spin_{T_w}(F_w) \xi_i^{-1}$.  

\begin{defn}
We say that a sequence of subsets $X_i$ of a topological space
is {\em becoming dense} if every open subset intersects $X_i$ for all sufficiently large $i$. 
\end{defn}

Put  $G = \spin_V(F_w), H = \spin_{T_w}(F_w)$.
We claim that if the closed subsets $\Gamma \backslash \Gamma \xi_i  H$ are becoming dense
in $\Gamma \backslash G$, then $W_i$ is spin globally representable for all sufficiently large $i$.  
Indeed, we need to check that $\Gamma \xi_i H \xi_i^{-1}$
intersects each $\K_w$ coset $\Gamma g \K_w$; note there are only finitely many possibilities
for $\Gamma g \K_w$.  Equivalently,
we need to check that 
that $\Gamma \xi_i H$ intersect $\Gamma g \K_w \xi_i$.
But, $\xi_i$ being constrained to a compact set, the number of possibilities
for $\K_w \xi_i$ is finite; so the latter statement would certainly follow if we know that $\Gamma \xi_i H$ are becoming dense.

\subsection{Application of the theorem of Ratner, Margulis/Tomanov; conclusion
of proof of Prop. \ref{prop:spin} and \ref{prop:spinbis}.} \label{Ratner}
As was indicated in the previous section, Prop. \ref{prop:spin}, or equivalently Prop. \ref{prop:spinbis},
will follow from the following statement: 

{\em Claim:} Let $V$ a quadratic space over $F$, $T_w$ a subspace of $ V \otimes F_w$ of dimension $\geq 7$, 
 $G = \spin_V(F_w)$, $H = \spin_{T_w}(F_w)$ and $\Gamma$ the arithmetic cocompact lattice in $G$
 defined in the previous section. 
Let $\xi_i \in G$ belong to a fixed compact subset of $G$ and have the property that
$\xi_i H \xi_i^{-1}$ is the stabilizer in $\spin_V(F_w)$ of a certain $F$-subspace $W_i$,
where no infinite subsequence of the $W_i$ have a nontrivial common intersection.
Then $\Gamma \backslash \Gamma \xi_i H$ is becoming dense in $\Gamma \backslash G$ 
as $i \to \infty$. 

%

%

%

%


To complete the proof of Prop. \ref{prop:spinbis}, we shall require two further results.
Firstly, we need a suitable consequence of the theorems of Ratner and Margulis-Tomanov; 
this is stated in Prop. \ref{prop:danimargulis} below and proved in the Appendix. 
The second is a group-theoretic
 result about generation of spin groups. 
  \begin{prop} \label{prop:danimargulis}
Suppose $\xi_i \in G$ remain within a compact set and, for any subsequence of $i$,
the subgroups $\xi_i H \xi_i^{-1}$ generate\footnote{It will suffice that they generate a Zariski-dense subgroup of $G$, as will be clear from the proof. } $G$.  Let $\mu_i$ be the $H$-invariant probability measure on $\Gamma \backslash \Gamma \xi_i H$. Then any weak limit of the measures $\mu_i$ is the $G$-invariant probability measure on $\Gamma \backslash G$. 
\end{prop}

\proof This is given in the Appendix. Note that one needs to verify the conditions
enumerated in Section \ref{conditions}; these follow from Lemma \ref{solvability} and standard facts. 
\qed 


 \begin{prop}\label{prop:grouptheory}
 Let $W_i$ be subspaces of a quadratic space $V$ over a nonarchimedean local field $F_w$
 and let $G = \spin_V(F_w)$. 
 Suppose that $\codim(W_i) \geq 7$ and that no infinite subsequence of the  $W_i$ have a common 
 nonzero intersection. Then
 the subgroup generated by the stabilizers of $W_i$ in $G$, is in fact all of $G$.  \end{prop}

The proof is given in Sec. \ref{sec:grouptheory}. 
Together these results prove immediately the {\em Claim} above, 
and therefore also  Prop. \ref{prop:spin},
\ref{prop:spinbis}. 

\subsection{Proof of Prop. \ref{prop:grouptheory}} \label{sec:grouptheory}

 In this section we give the proof of Prop. \ref{prop:grouptheory}. 
During this section, we will work exclusively with the $F_w$-points of certain algebraic groups over $F_w$.
Consequently, for brevity, we write simply (e.g.) $O_V$ or $\spin_V$ instead of
$O_V(F_w)$ or $\spin_V(F_w)$.

We will need a few lemmas.  As before we set $T_i = W_i^{\perp}$;
then no infinite subsequence of the $T_i$ are contained in a common
proper subspace.  Recall that by $O_{T_i}$ we mean the stabilizer of
$W_i$ in the orthogonal group $O_V$.  We will first prove that the subgroup generated by $O_{T_i}$ is all of $O_V$, and finesse
the claimed result from this.  
Let $\Xi$ be the subgroup generated by all the $O_{T_i}$.

We shall first check that $\Xi$ acts transitively on vectors in $V$
of length $1$.

Let $U$ be any nondegenerate subspace of $V$ of dimension at least
$7$, and $\proj$ the orthogonal projection from $V$ 
to $U^{\perp}$. Recall that $r_w$ denotes the reflection through the hyperplane
perpendicular to $w$, whenever $w \in V$ is such that $q(w) \neq 0$. 

\begin{lemma}
Let $v \in V$ and $w \in V$ so that $q(w) \neq 0$. 
Suppose $w \notin U \cup U^{\perp}$ and $v \notin U^{\perp}$. There exists $g \in \langle O_U, r_w \rangle$
such that $\proj(gv) - \proj(v) = \proj(w)$ and $g v \notin U^{\perp}$. 
\end{lemma}
\proof
We will find such a $g$ of the form $r_w \sigma$, for an appropriate choice
of $\sigma \in O_U$ to be made at the end.   
Write $v = u_0 + u_0^{\perp}$, with $u_0 \in U, u_0^{\perp} \in U^{\perp}$. 
Note that $gv \notin U^{\perp}$ if $\sigma v \notin \langle U^{\perp}, w\rangle$, 
which is certainly true if $\sigma u_0$ does not belong to the line spanned
by the projection of $w$ to $U$.

Now
$$ \proj( g v - v) = \proj(r_w \sigma u_0 ) + \proj(r_w  u_0^{\perp} - u_0^{\perp})$$
Now $r_w u_0^{\perp} - u_0^{\perp}$ is a certain multiple of $w$. 
We therefore want to solve $\proj(r_w \sigma u_0 ) = \proj(w) -  \proj(r_w  u_0^{\perp} - u_0^{\perp})$;
the right hand side is certainly a multiple of $\proj(w)$. Moreover, we want $\sigma u_0$
not to belong to the line spanned by the projection of $w$ to $U$.

The map $u \mapsto \proj(r_w u)$ is, by assumption that $w \notin U \cup U^{\perp}$, a surjection from $U$ onto the line
spanned by $\proj(w)$. Let $K$ be its kernel. 

$K$ is a certain subspace of $U$ of codimension $1$. Choose a nondegenerate subspace $J \subset K$
of codimension at most $2$ inside $U$. (Let $K^{\perp}$ be the orthogonal complement to $K$ within $U$. If $K^{\perp} \cap K$ is trivial, we can take $J = K$;
otherwise, $K^{\perp} \subset K$, and one can take for $J$ any codimension $1$ subspace of $K$
not containing $K^{\perp}$). 

We claim that every $J$-coset (and consequently every $K$-coset) has
nonempty intersection with every level set of the form $\{u \in U :
q(u) = c\}$ Indeed, it suffices to check that for any linear
functional $l$ on $J$ and for any $c' \in F_w$, the equation $q(u) +
l(u) = c'$ is solvable with $u \in J$; but because $J$ is
nondegenerate, we can convert this (after an affine change of
coordinates) to the equation $q(u) = c''$. This is solvable as long as
$\dim(J) \geq 5$, by Lemma \ref{solvability}.  More precisely, the
intersection of the $J$-coset with the level set is the $F_w$-points
of an affine quadric of dimension at least $4$.  But such a quadric is
automatically isomorphic over $F_w$ to an open subset of projective
space; in particular, its $F_w$-points are Zariski dense.

It follows that there exists $u \in U - \{0\}$ satisfying 
$q(u) = q(u_0)$ and $\proj(r_w u) = \proj(w) - \proj(r_w  u_0^{\perp}
- u_0^{\perp})$. Moreover (by the Zariski density) 
this $u$ can be chosen to avoid the line spanned by the projection of $w$ to $U$.
Now choose $\sigma$ so that $\sigma u_0 = u$. 
\qed 

\begin{lemma} \label{transitiveunit}
$\Xi$ acts transitively on vectors $v \in V$ satisfying $q(v) = 1$. 
\end{lemma}

\proof

Take $v_1, v_2$ with $q(v_1) = q(v_2) = 1$.
By virtue of the assumption that no infinite subsequence of $T_i$ are contained
in a proper subspace, there exists $T_i$ perpendicular to neither $v_1$ or $v_2$; call it $U_0$. 
Let $\proj$ be the orthogonal projection onto $U_0^{\perp}$. 

Again, because no infinite subsequence of $T_i$ are contained in a proper subspace, the subspaces $T_i$ not contained in $U_0$ or $U_0^{\perp}$ 
 span $V$. Therefore, the subspace of $U_0^{\perp}$ spanned by $\proj(t_i)$ with $t_i \in T_i, q(t_i) \neq 0$
 and $t_i \notin U_0 \cup U_0^{\perp}$
is a linear subspace of $U_0^{\perp}$ that is also topologically dense; so it coincides with $U_0^{\perp}$. 

It follows that there exist a finite list $t_i \in T_i$ with $q(t_i) \neq 0, t_i \notin U_0 \cup U_0^{\perp}$ for all $i$ and so that 
$\proj(v_1) - \proj(v_2) = \sum_{i} \proj(t_i)$.  Repeated applications of the previous Lemma show
that there is $g \in \Xi$ with $\proj(g v_1) = \proj(v_2)$ and $g v_1 \notin U_0^{\perp}$. 
Thus the projection of $g v_1$ and $v_2$ to $U_0$ both have the same norm and are both nonzero. 
Modifying $g$ by an element of $O_{U_0}$ as necessary, 
we see that $\Xi$ maps $v_1$ to $v_2$ as required. 
\qed

Now, $\Xi$ by definition is the subgroup generated by the $O_{T_i}$.
It follows from this and Lemma \ref{transitiveunit} that $\Xi$ contains
all reflections through vectors of length $1$ (for there exists such a vector in any of the $T_i$).  By Lemma \ref{solvability},
$\Xi$ contains the image of $\spin_V$ in $O_V$.  Since (by Lemma \ref{solvability} again)
the spinor norm is surjective on $\so_{T_i}$, we deduce from this that $\Xi$ contains
also $\so_V$, and so $\Xi=O_V$. 

We are now ready to complete the proof of Prop. \ref{prop:grouptheory}. 

\proof
Let $\mathscr{P}$ be the class of hyperbolic planes inside $V$. 
For any $P \in \mathscr{P}$ and any $T_i$, we first claim that the orbit
$O_{T_i} \cdot P \subset \mathscr{P}$ coincides with the orbit $\spin_{T_i} \cdot P$. 
For this, it suffices to check that if $H_P$ is the stabilizer of
$P$ in $O_{T_i}$, then $H_P$ surjects onto the finite quotient
$O_{T_i}/\spin_{T_i}$. But $H_P$ contains the pointwise stabilizer, in $O_{T_i}$,
of the projection of $P$ to $T_i$; it is easy to see that $H_P$
 contains the orthogonal
group of a nondegenerate subspace $W_0$ of $T_i$ of codimension $\leq 4$.
Since $\dim(T_i) \geq 7$, we see that the dimension of $W_0$ is $\geq 3$;
then $O_{W_0}$ contains a reflection and the spinor norm is surjective on $SO_{W_0}$
by Lemma \ref{solvability}. 

Let $\Xi_0$ be the subgroup of $\spin_V$ generated by the stabilizers of the $W_i$. Then it follows from the remark above that the $\Xi_0$-orbits on 
$\mathscr{P}$ coincide with the $\Xi$-orbits. But, by Witt's theorem, 
$\Xi = \ortho_V$ acts transitively on $\mathscr{P}$.  So $\Xi_0$ also acts transitively on $\mathscr{P}$. 

Now there exists at least one $P \in \mathscr{P}$ that is a subspace of some $T_i$, because,
since $\dim(T_i) \geq 7$, $T_i$ is isotropic (Lemma \ref{solvability}). 
So $\Xi_0$ contains $\spin_P$ for this choice of $P$.
Thus $\Xi_0$ contains the subgroup of $\spin_V$ generated by $\{\spin_P: P \in \mathscr{P}\}$. 
This  coincides with $\spin_V$ by Lemma \ref{solvability}. 
\qed

\section{Algebraic structures associated to integral orbits.}\label{algebraicstructures}
This section is an expansion of the brief remarks in
Section \ref{stabilizer}. It is devoted to a discussion of ``class number problems,''  and the role
of the stabilizer. We have included this material since we believe it gives the correct context
for our work; however, we note that the proof of Theorem~\ref{th:main}
is independent of the material in this section. 

Much of the material in this section may be found, implicitly or explicitly and in various contexts,
in the work of many authors (see Section \ref{connection} for some references). Indeed, the study of class number problems begins with Gauss. Our goal 
in this section has been to give a coherent discussion of such problems in a quite general setting.

\subsection{Class number problems.} \label{classnumber}

Let $\mathbf{G}$ be an algebraic group over $\Q$; let $\rho:
\mathbf{G} \to \GL(V)$ be a $\Q$-linear representation of
$\mathbf{G}$.  Let $v_0 \in V$ be so that the orbit $\mathbf{G} . v_0$
is a closed variety $\mathbf{X}$; let $\mathbf{H}$ be the isotropy
subgroup of $v_0$. $\mathbf{H}$ must be reductive; indeed, this is equivalent
to $\mathbf{X}$ being affine.   Let $\Gamma$ be a
congruence lattice in $G := \mathbf{G}(\mathbb{R})$ and let $\Lambda$ be a
$\Gamma$-stable lattice in $V$. We set\footnote{We note that $X_{\Z}$ may be empty without the following discusion becoming vacuous; indeed, 
the following discussion may be used to {\em prove} $X_{\Z}$ is nonempty, as is done in the text.} $X_{\Z} := \mathbf{X}(\Q) \cap
\Lambda$.


  It is worth remarking that all the morphisms of groups we
discuss depend on the choice of $v_0 \in X$.\footnote{It would be
instructive to understand how the objects in this paper vary with
choice of basepoint, and whether a more canonical construction is
possible.  For instance: the solutions to $x^2 + y^2 + z^2 = d$, as we
have seen, are in bijection with a class group, but not canonically
so; the canonical structure on this set is that of a torsor for a
class group.}

{\em Class number problem:} Understand the parameterization and distribution
of $\Gamma$-orbits on $X_{\Z}$.

In Section \ref{param1} we carry this out in a rather {\em ad hoc} way, first parameterizing
orbits over $\Q$ and then passing to $\Z$. In Section \ref{param2} we describe
a more unified, but less concrete approach in terms of torsors.  Either way, the material of this Section is intended to justify
the approximate discussion of Section \ref{stabilizer}. Although to apply
either of the parameterizations (of Section \ref{param1} or Section \ref{param2})
requires (possibly complicated) local computations, 
the discussion still has considerable explanatory power at a heuristic level. 
For example, the situations where the $\Gamma$-orbits on
$\mathbf{X}(\Z)$ can be given the structure of a group (or at least a
principal homogeneous space for a group) should be precisely those
where the stabilizer $\G_x$ is abelian, as in the case of
representations of integers by ternary forms discussed in Section \ref{linnik}. 

\subsection{Parametrization of orbits} \label{param1}
 Let $\adele_f$ be the ring of finite adeles
and let $K_f$ be the closure of $\Gamma$ in $\G(\adele_f)$; thus
$\Gamma = \G(\Q) \cap K_f$.   The quotient $K_f \backslash
\G(\adele_f)/ \G(\Q)$ is finite;
we refer to it as the {\em genus} of $\Gamma$. 
Fix a set of representatives $\{1=g_1, \dots, g_h\}$ for the cosets. For each $1 \leq i \leq h$, 
set $\Lambda_i = g_i^{-1} \cdot \Lambda$ (recall that $\G(\adele_f)$ acts naturally on lattices in $V$)
and $\Gamma_i = \G(\Q) \cap g_i^{-1} K_f g_i$. Then $\Lambda_i$ is stable under $\Gamma_i$,
and $\Lambda_1 = \Lambda, \Gamma_1 = \Gamma$. We will eventually describe
a parameterization of $\bigcup_{i} \Gamma_i \backslash (\mathbf{X}(\Q)
\cap \Lambda_i)$.  As is common in this genre of problem, it is
significantly easier to understand this union of orbits than the
individual orbits themselves.

\subsubsection{Parametrization of $\Q$-orbits.} The orbits of $\G(\Q)$ on $\mathbf{X}(\Q)$ are parameterized
by the kernel of the map of pointed sets $$H^1(\Galois, \mathbf{H}(\overline{\Q})) 
\to H^1(\Galois, \mathbf{G}(\overline{\Q})).$$ 
Explicitly, given a representative $x \in \mathbf{X}(\Q)$ for an orbit, there exists $g \in \G(\overline{\Q})$
such that $g v_0 = x$; for each $\sigma \in \Galois$, we have $g^{\sigma} v_0 = x$, 
so that $\sigma \mapsto g^{-1} g^{\sigma}$ defines an element of the former set.  
\subsubsection{Parametrization of $\Z$-orbits within a $\Q$-orbit.} 
We fix a $\Q$-orbit $\G(\Q) . x$ and are interested in parameterizing
$\Gamma$-orbits on $\G(\Q) .x \cap X_{\Z}$. Let $\class$ be the set of such classes. 
More generally, let $\class_i$ be the set of $\Gamma_i$-orbits on $\G(\Q). x \cap \Lambda_i$.
Thus $\class_1 = \class$. As is usual in such problems, it will be simpler to parameterize
$\bigcup_{i} \class_i$, the union of classes associated to the genus of $\Gamma$.

Let $\G(\adele_f)^{\integer} = \{g \in \G(\adele_f): g . x \in \Lambda \otimes \hat{\Z}\}$;
thus $\G(\adele_f)^{\integer} = \prod_{p} \G(\Q_p)^{\integer}$, where 
$\G(\Q_p)^{\integer} :=
\{g_p \in \G(\Q_p): g . x \in \Lambda \otimes \Z_p\}$.
We note that $K_f . \G(\adele_f)^{\integer} = \G(\adele_f)^{\integer}$. 
Then the set of $\delta \in \G(\Q)$ such that $\delta . x \in X_{\Z}$
is exactly $\G(\Q) \cap \G(\adele_f)^{\integer}$. It follows that $\class$ is naturally identified with the quotient
$\G(\Q) \cap K_f \backslash \G(\Q) \cap \G(\adele_f)^{\integer} / \G_x(\Q)$, that is to say: 
$$\class \stackrel{\sim}{\rightarrow} K_f \backslash K_f \G(\Q) \cap \G(\adele_f)^{\integer}/ \G_x(\Q)$$
Similarly, $$\class_i \stackrel{\sim}{\rightarrow} K_f \backslash K_f g_i \G(\Q) \cap \G(\adele_f)^{\integer}/\G_x(\Q).$$ We conclude that the union $\bigcup_{i} \class_i$ is naturally identified with
$$\bigcup_{i} \class_i \stackrel{\cong}{\to} K_f \backslash \G(\adele_f)^{\integer}/ \G_x(\Q)$$
Thus $\bigcup_{i} \class_i$ is easy to describe: it maps to $K_f \backslash \G(\adele_f)^{\integer} /\G_x(\adele_f)$, 
which can be computed purely locally; and the fiber above the class $K_f . q. \G_x(\adele_f)$
is identified with $q^{-1} K_f q \cap \G_x(\adele_f) \backslash \G_x(\adele_f) / \G_x(\Q)$:

$$\bigcup_{i} \class_i = \bigcup_{ q \in K_f \backslash \G(\adele_f)^{\integer}/\G_x(\adele_f) }
q^{-1} K_f q \cap \G_x(\adele_f) \backslash \G_x(\adele_f) / \G_x(\Q).$$

Note that $K_f \backslash \G(\adele_f)^{\integer}/\G_x(\adele_f)$
is precisely the set of $K_f$ orbits on $\G(\adele_f). x \cap (\Lambda \otimes \hat{\Z})$.  The computation of this set of orbits is a purely {\em local} problem.

Thus what we have shown may be phrased: the union of classes $\class_i$,
where $i$ varies over the genus of $\Gamma$, is parameterized by 
a certain union of adelic quotient spaces associated to $\G_x$. 

 \subsubsection{A diagram.}
We can summarize this discussion in the following diagram, where the left-hand
vertical sequence of sets is exact in the sense that the first term is
exactly the fiber over an element $q$ of the last term.

\begin{equation}
\begin{CD}
q^{-1} K_f q \cap \G_x(\adele) \backslash \G_x(\adele) / \G_x(\Q) & @>g \mapsto qg>> 
& K_f \backslash \G(\adele_f) /\G(\Q) \\
 @VVV & & @V=VV  \\ 
\bigcup_{i} \class_i & @>>> & i \in \{1,2,\dots, h\} \\
@VVV & &   \\ 
\{ q \in K_f \backslash \G(\adele_f)^{\integer} /\G_x(\adele_f)\} 
&  & \\
\end{CD}
\end{equation}

In practice, this subdivides the study of the original orbit set $\class = \class_1$ into two sub-problems: 
\begin{enumerate}
\item Local problem: understand the set of orbits $K_f \backslash \G(\adele_f)^{\integer}/\G_x(\adele_f)$. 
\item Global problem: in order to ``recover'' $\class$ from the union $\cup_{i} \class_i$, we must understand the behavior of the maps from 
$q^{-1} K_f q \cap \G_x(\adele) \backslash \G_x(\adele) / \G_x(\Q) $ to 
$K_f \backslash \G(\adele_f) /\G(\Q) $. It is this which can be approached via ergodic methods, 
for it is evidently related to the dynamics of the action of $\G_x(\adele)$ on $\G(\adele_f)/\G(\Q)$. 
\end{enumerate}

We note that the term at the bottom left is trivial if $K_f$ acts
trivially on $\X(\Q_p) \cap (\Lambda \otimes \Z_p)$ for each $p$, i.e.,
there is {\em locally} only one orbit on integral points. 
  In this case,
if we can show that the top horizontal map of adelic quotients is
surjective, we will have shown that $\mathbf{X}(\Q) \cap \Lambda_i$ is
nonempty for all $i$.  As we shall see, in the context of
representations of quadratic forms this will show exactly that a form
$Q'$ is represented by, not only {\em some} form in the genus of $Q$,
but {\em every} form in the genus of $Q$.



\subsection{Examples}

 We give several examples
but do not carry out the local computations in any detail. 

\begin{enumerate}
\item Quadratic forms. Let $Q$ be a quadratic form on
the $\Z$-lattice $\Lambda$ of rank $n$. Let $V = \Lambda \otimes \Q$,  $\G = \so(Q)$ and $\Gamma$
the stabilizer of $\Lambda$ in $\G(\Q)$; $\rho$ is the natural representation of $\G$ on $V$.

Let $0 \neq d \in \Z$. The level set $Q(\mathbf{x}) = d$ is a closed subvariety $\X \subset V$
which is a homogeneous space for $\G$. 
Witt's theorem shows that $\G(\Q) $ acts transitively on $\X(\Q)$. 
The stabilizer $\G_x$ of a point $x \in \X(\Q)$ is the orthogonal group $\o(\langle x \rangle^{\perp})$,
and the adelic quotient $q^{-1} K_f q \cap \G_x(\adele) \backslash \G_x(\adele) / \G_x(\Q)$ is closely related
to the genus of the quadratic form induced on $\langle x \rangle^{\perp}$.

The considerations of the previous section show that {\em the
  $\Gamma$-orbits on representations of $d$ by a quadratic form in $n$
  variables} are closely related to the {\em genus of a certain
  collection of quadratic forms in $n-1$ variables.} 

This observation is, of course, not new and seems to be classical. It is quite explicitly presented in Kneser's article \cite{Kneser}. Shimura's book \cite{shimura} carries out some of the difficult local computations associated
to precisely implementing this. 

Some particular and familiar corollaries of this observation are:
\begin{enumerate}
\item The case previously discussed, and due to Gauss: that the number of representations
of $n$ by the form $x^2+y^2+z^2$ is related to the class number of $\Q(\sqrt{-n})$
(and therefore to genera of quadratic forms of rank $2$). 
\item If the signature of $Q$ is $(n-1,1)$, and $n \geq 4$, then the representation numbers show the following curious behaviour: the number of orbits on $Q(\mathbf{x}) = d$, as $d$ varies
through squarefree integers, grows roughly as $|d|^{\frac{n}{2} - 1 \pm \varepsilon}$
as $d \to -\infty$; on the other hand it
grows very slowly (say as $|d|^{\varepsilon}$) as $d \to \infty$.
The difference is that the stabilizer $\G_x$ in the former case is compact at $\infty$,
and in the latter case is semisimple and noncompact at $\infty$, therefore satisfying strong approximation. \footnote{Recall that if $\H$ is a semisimple, simply connected $\Q$-group
for which $\H(\mathbb{R})$ is noncompact, then $\H(\Q)$ is dense in $\H(\adele_f)$, and in particular
the quotient space $\Omega \backslash \H(\adele_f)/\H(\Q)$ is a singleton for any open
subgroup $\Omega$. Even if $\H$
fails to be simply connected, this set is parameterized in a completely understandable way
by the center $\mathbf{Z}(\H)$ and in particular has ``very few'' elements.}
\end{enumerate}

\item Class groups of number fields. Let  $V$ be an $n$-dimensional $\Q$-vector space together
with a $\Z$-lattice $V_{\Z}$, $\G = \SL(V)$ and $\Gamma = \SL(V_{\Z})$, 
and consider the representation $\rho$ of $\G$ on the vector space $W = \mathrm{Sym}^n(V^{*})$ of homogeneous polynomials
of degree $n$ on $V$, defined by
$\rho(g) f = f(x g) $. 
Then $\Lambda = \mathrm{Sym}^n(V_{\Z}^*)$, a lattice in $W$ that is preserved by $\rho(\Gamma)$. 

Suppose $w \in W$ is a degree-$n$ form which factors over $\bar{\Q}$ into a
product $\ell_1 \cdot \ldots \cdot \ell_n$ of linear forms; we refer to the
square of the determinant of the resulting element of $\Hom(V,\Z^n)$ as the {\em
  discriminant} of $w$.  For each nonzero integer $d$, write $\X_d$ for
the closed subvariety of $W$ parametrizing forms which factor into linear
forms over $\bar{\Q}$ and have discriminant $d$.  Then $\X_d$ is a
homogenous space for $G$, and the stabilizer in $G$ of a point in $\X_d$ is
the semidirect product of a torus with a finite group scheme
geometrically isomorphic to $A_n$.  We call a form in $\X_d$ {\em
  primitive} if, for every prime $p|d$, the reductions of $\ell_1, \ldots,
\ell_n$ have the property that every subset of size $n-1$ intersects
transversely.  Then the $\Gamma$-orbits of primitive points of $\X_d(\Z)$
will parametrize classes in Picard groups of certain orders of degree $n$ over
$\Z$ with discriminant $d$.  When $n=2$, this reduces to the classical
Gauss correspondence between binary quadratic forms and narrow class groups
of quadratic orders.  When $n=3$, this is quite close to the case considered in the
introduction to \cite{Bhargava3}, which eventually replaces this space
with a simpler one and thereby provides a concrete parametrization of
ideal classes in cubic rings.


\item M. Bhargava in his sequence of papers (\cite{Bhargava}, \cite{Bhargava2}, \cite{Bhargava3}) has studied many lovely examples
of class number problems in the case where the $\G$-action is {\em prehomogeneous}: that is to say, 
the ring of invariants for $\G$ acting on $V$ is generated by a single
polynomial $f$.  In these cases, Bhargava develops new composition
laws and relates the classification of integral orbits to various
structures in algebraic number theory (e.g. class groups of orders and
$n$-torsion in class groups of orders.)  For example, he studies the
action of $(SL_2)^3$ on the space $(\Z^2)^{\tensor 3}$ and shows that the
orbits are related to triples of ideal classes for quadratic orders
with product $1$.  A remarkable feature of his constructions is that
they completely deal with the (rather complicated) local problems
implicit in our discussion above; it is to avoid these local problems
that we have restricted our attention to representations of forms with
squarefree discriminant in the present paper.
 \end{enumerate}

\subsection{Orbits over more general bases: relation with torsors} \label{param2}
It is an interesting open problem to understand the extent to which
the framework of ``class group problems'' can be generalized to base
schemes other than $\Spec \,\Z$.  Certainly it is well-known that some
version of Gauss composition for quadratic forms can be carried out
over an arbitrary commutative ring (see, e.g., Kneser
\cite{knes:jnt}).  Because we will not need to work over an arbitrary
base in the present paper, we will confine ourselves to a few speculative remarks
here.
 
One possibility for the general set-up is as follows.  Let $S$ be a scheme, $X \ra S$ an fppf morphism, and $G \ra S$
an fppf group scheme.  Suppose that $G$ acts on $X$; this action
defines and is defined by a morphism
\beq
m: G \times X \ra X \times X
\eeq
defined by $m(g,x) = (gx,x)$.  Suppose that $m$ is also fppf; in this
case, we say the action of $G$ on $X$ is faithfully flat.   

Now choose a basepoint $x_0 \in X(S)$, and let $H \subset G$ be the stabilizer
of $x_0$.  Then for any other $x \in X(S)$ one can define the space of
paths $P_{x_0,x}$ to be $m^{-1}(x,x_0)$. The association $x \mapsto P_{x_0,x}$ assigns to every $x$ an fppf $H$-torsor over $S$; evidently, if $x$ and $y$ are in the same $G(S)$-orbit, the $H$-torsors
$P_{x_0, x}$ and $P_{x_0,y}$ are isomorphic. 
So one gets a map from the set of orbits $G(S) \bs
X(S)$ to the fppf cohomology group $H^1(S,H)$, whose image is just the
kernel of the natural map $H^1(S,H) \ra H^1(S,G)$.  In particular, if
$H^1(S,G)$ is trivial and $H$ is abelian, the orbit set acquires the
structure of a group. This is very likely related to the composition laws presented in 
\cite{Bhargava2}, \cite{Bhargava3}, and seems likely to suggest
further composition laws on integral orbit spaces.

There are several potential advantages to studying the problem of
integral orbits in this generality.  For instance, in the case $S = \Spec \, \Z$, 
\begin{itemize}
\item The cohomology set $H^1(\Spec \, \Z, H)$ incorporates, in one step,
  the Galois-cohomological data recorded by $H^1(\Spec \, \Q, H)$, and
  the adelic data recorded by the kernel of $H^1(\Spec \, \Z, H) \ra
  H^1(\Spec \, \Q, H)$.  For instance, in the case treated by Bhargava
  in \cite{Bhargava3}, where $G = GL_2 \times GL_3 \times GL_3$ and $X
  = \Z^2 \tensor \Z^3 \tensor \Z^3$, the class in $H^1(\Spec \, \Q,
  H)$ keeps track of a cubic field, while the extra data coming from
  $H^1(\Spec \, \Z, H)$ yields an ideal class in some order of that field.

\item Certain restrictions are imposed on us by the requirement that
  the multiplication map $m$ be flat.  This condition implies in
  particular that the stabilizer of any point $x \in X(S)$ is flat over
  $S$.  
  
  For instance, if $X_d \subset \A^3$ is the space of binary quadratic forms
  of discriminant $d$, then the action of $\SL_2$ on $X_d$ need not be
  flat; if $p^2 | d$, and $x$ is a form in $X_d(\Z)$ which reduces to
  $0$ mod $p$, then the stabilizer of $x$ is evidently not flat.  One
  can fix this problem by considering instead the quasi-affine scheme
  obtained by removing the origin from $X_d$.  In more classical
  language, we have restricted our attention to primitive quadratic
  forms.  In general, a natural candidate for the correct notion of  ``$x \in X(S)$ is primitive'' for an integral orbit problem over a
  base $S$ should be ``the stabilizer of $x$ is flat over $S$.''  One
  nice feature of Bhargava's work is that it does not restrict itself
  to primitive situations.  As one might expect from the above
  discussion, Bhargava typically finds a subset of primitive orbits
  among the set of all orbits which admit a composition law, while the
  set of all orbits does not.  (For example, in \cite{Bhargava3}, the
  composition law on $2 \times 3 \times 3$ cubes applies precisely to
  those cubes which are {\em projective} in Bhargava's sense.)

\end{itemize}

We remark, finally, that this framework is natural for
understanding the manner in which various classical constructions
depend on choice of basepoint $x_0$; rather than fixing a basepoint,
it is probably best to consider the gerbe $G \backslash X$; any choice
of $x_0 \in X(S)$ provides an isomorphism between this gerbe and the classifying stack of
the stabilizer of $H_{x_0}$, but there is no such canonical
isomorphism in general.

\section{Extensions and problems}
As remarked, the methods used to prove Theorem \ref{thm:final} can be extended and optimized in several ways.

 It is applicable
also to other embedding problems (e.g., pertaining to hermitian forms; a slightly more ``exceptional'' example is the embeddings of a cubic
order into matrix algebras over the octonions, considered in \cite{GG}) as well as to other equidistribution problems (for instance, one can expect to understand, by this technique, the distribution of all integral, positive definite, quadratic forms of discriminant $D \rightarrow \infty$ inside the moduli space $\PGL_n(\Z) \backslash \PGL_n(\mathbb{R}) / PO_n$
of homothety classes of quadratic forms; the case $n=2$ is a theorem of W. Duke, whereas
for $n > 2$ and {\em indefinite} quadratic forms, the analogous result is due to A. Eskin and H. Oh; the $p$-adic methods of this paper allow the treatment of the outstanding case). As remarked previously, we also have not optimized the results even for quadratic forms; the condition $n \geq m+7$ is not the limit of the method, and our method should also yield an asymptotic for the representation numbers.

Let us remark on some more amibitious extensions and problems:
\begin{enumerate}
\item Effectivity; bounds for Fourier coefficients of Siegel modular forms.

As remarked, a fundamental defect of Theorem \ref{th:main} is its ineffectivity. This arises from the ineffectivity of Ratner's theorem. Were one to have, in the context of Proposition \ref{prop:danimargulis}, an effective estimate
on the rate of convegence of the $\mu_i$ to their limit, this would yield an effective version of Theorem \ref{th:main}.  While it is plausible that existing proofs of Ratner's theorem may be effectivized, 
a much bigger challenge is to obtain a reasonable rate of convergence (say, polynomial in the relevant parameters). 

In this context, it should be noted that Margulis has given a beautiful effective proof of the convergence
of the invariant measures on closed $\mathrm{SO}(2,1)(\mathbb{R})$-orbits on $\SL_3(\Z) \backslash \SL_3(\mathbb{R})$ 
to their limit. Although the present situation is quite different, and more complicated
(because there are many intermediate subgroups) this result certainly makes it plausible that an effective result is possible.

In any case, another significant payoff of such an effective result would be a nontrivial estimate on the Fourier coefficients of Siegel modular forms arising from $\theta$-series of quadratic forms.

\item Representations in codimension $2$.

We have remarked at various points that the natural limit of the method presented
in this paper is $n  = m+3$; for, in the case $n = m+2$, one is forced
to consider actions of a $p$-adic {\em torus} $\mathrm{SO}_2(\Q_p)$
on a homogeneous space $\Gamma \backslash \SO_n(\Q_p)$.  

However, there is nevertheless
a natural approach to the case $n=m+2$, replacing our use of Ratner theory
by the emerging theory of torus rigidity (see, in particular, the survey and announcement
\cite{EL}). The idea is to replace the use of $\Q_p$ by a product of two completions
$\Q_{p} \times \Q_q$, and consider the action of $\mathrm{SO}_2(\Q_p) \times \mathrm{SO}_2(\Q_q)$
on $\Gamma \backslash \SO_n(\Q_p) \times \SO_n(\Q_q)$. For suitable $p$ and $q$, 
the group $\mathrm{SO}_2(\Q_p) \times \mathrm{SO}_2(\Q_q) \cong \Q_p^* \times \Q_q^*$
is a ``higher rank'' torus, and one expects a certain degree of ridigity for the invariant ergodic measures. 

There are several obstacles to this approach. For one, the relevant measure rigidity statements
are not (yet) available. Moreover, they require a pre-condition: {\em positive entropy}. Another more serious obstacle is that, in the torus case, one does not have
a good way of ruling out concentration of limit measures on intermediate subgroups.

Nevertheless, it does not seem entirely impossible that these obstacles can be overcome. We refer, in particular,
to the series of papers \cite{ELMV} where essentially the analogous question is considered, but replacing $\SO_n$ with $\PGL_n$ and $\SO_2$ with a maximal torus, and it is shown how to overcome
these obstacles in several situations. In particular, satisfactory results are obtained for $n=3$. 
\end{enumerate}

 \appendix{}
 \section{}
 We now give the proof of Proposition \ref{prop:danimargulis}.  
 
 The ideas follow the ``linearization''
 technique which we learned from \cite{DM}; however, we simply the computations considerably by 
 Lemma \ref{ergodiclimit}. 
This Lemma was noted, in an entirely different context, by the second author jointly with M. Einsiedler and E. Lindenstrauss. It was pointed out to us by Y. Shalom that it appears already in the paper of Glasner and Weiss \cite{GW}. 
 
 It is important to note that, while the ``trick'' of Lem. \ref{ergodiclimit} makes the proofs much easier,
 the original ideas of \cite{DM} carry over to the $p$-adic setting without essential change, and this would be needed to treat the case where $H$ does not have property (T) (notation of Prop. \ref{prop:converge}). 
 \subsection{Ergodicity of limit measures for a group with $T$.}

 \begin{lemma}  \label{ergodiclimit} \cite{GW}  Let $H$ be a locally compact, second countable
 group with property $(T)$. Let $\mu_i$ be a sequence of ergodic $H$-invariant probability measures on a locally compact space $X$. Then any weak limit  $\mu_{\infty}$ of the $\mu_i$ is also an ergodic $H$-invariant measure. 
 \end{lemma}
 ( We note that property $(T)$ could be replaced for a uniform spectral gap for the $H$-action on the representations $L^2(X, \mu_i)$. This is not an entirely idle comment, as it allows one to apply
 the same reasoning in many cases when $H$ has rank one, by suitable bounds on automorphic spectrum.) 
 
 \proof 
In the interest of self-containedness, we present a proof, at least in the case when $H$ is a {\em discrete} group. This suffices for the present application (for we apply it when $H$ is a $p$-adic Lie group
which admits a lattice with property (T)); however, the proof is easily modified to handle the general case. 
 
Thus, take $S$ to be a generating set for $H$. By Property (T), we may choose $\delta > 0$ so that for any unitary $H$-representation $\rho \rightarrow \mathrm{Isom}(V)$ on a Hilbert space $V$ not containing the trivial representation, and for any $0 \neq v \in V$, we have
$$\sup_{s \in S} \|\rho(s) v - v\|  > \delta \|v\|$$
Let $T := \frac{1}{2|S|}\sum_{s \in S} (s+s^{-1}) \in \C[H]$, the group algebra of $H$. It follows that there is some $\beta < 1$ -- depending only on $\delta$ -- such that $\|T v\| \leq \beta \|v\|$ for all $v \in V$, where $V$ is as above.

Let $f, g \in C_c(X)$, the space of continuous compactly supported functions on $X$, and write $\bar{f}_i = f - \int_{X} f(x) d\mu_i(x)$. 
We note that $\|\bar{f}_i\|_{L^2(\mu_i)} \leq \|f\|_{L^2(\mu_i)}$. Clearly
$$ \langle T^n f, g \rangle_{L^2(\mu_i)} = \langle T^n \bar{f}_i, \bar{g}_i \rangle _{L^2(\mu_i)}
+ \int f d\mu_i \int g d\mu_i.$$
Note that, using the ergodicity of $\mu_i$, we have the bound $$\langle T^n \bar{f}_i, \bar{g}_i \rangle_{L^2(\mu_i)} \leq \beta^n \|f\|_{L^2(\mu_i)} \|g\|_{L^2(\mu_i)}.$$ Therefore 
\begin{equation} \label{mix} \left| \int T^n f(x) g(x) d\mu_i(x) - \int f(x) d\mu_i \int  g(x) d\mu_i \right|  \leq  \beta^n \|f\|_{L^2(\mu_i)} \|g\|_{L^2(\mu_i)}\end{equation} 

The assertion (\ref{mix}) passes to the limit $i \rightarrow \infty$, and the corresponding assertion holds also replacing
$\mu_{i}$ by $\mu_{\infty}$.  It then extends from $C_c(X)$ to $L^2(X)$ by density. 
Thus for $f,g \in L^2(X, \mu_{\infty})$ we have:
\begin{equation} \label{mix2} \left| \int T^n f(x) g(x) d\mu_{\infty}(x) - \int f(x) d\mu_{\infty} \int  g(x) d\mu_\infty \right|  \leq  \beta^n \|f\|_{L^2(\mu_\infty)} \|g\|_{L^2(\mu_\infty)}\end{equation} 
If $S$ is a $H$-invariant measurable subset for $\mu_{\infty}$, take $f=g=1_S$, the characteristic function,
and take $n \rightarrow \infty$ to see that $\mu_{\infty}(S)^2 = \mu_{\infty}(S)$, i.e. $S$ is null or conull. 
\qed

\subsection{Growth properties of $p$-adic polynomials}
Let $F$ be a number field, $w$ a nonarchimedean place of $F$, and $F_w$ the corresponding completion.

 For any $k > 0$ and $\mathbf{x} = (x_1, \dots, x_k) \in F_w^k$ we put
 $\|\mathbf{x}\| = \sup_{i} |x_i|$.  A {\em ball} in $F_w^k$ is
 a subset of the form $\{\mathbf{x}: \|\mathbf{x} - \mathbf{x}_0\| \leq \delta\}$. 
 The ultrametric property assures us that two balls are disjoint or one is contained in the other.  
 Put $\order_w[M] = \varpi_w^{M} \order_w$, where $\varpi_w$ is a uniformizer at $w$.

 \begin{lemma}
 Let $t > 0$. 
 Let $\theta: \order_w \rightarrow \order_w^k$ be a polynomial map of degree $d$
so that $\sup_{\lambda \in \order_w} \|\theta(\lambda)\|  \geq t $. 
 Then there is a continuous function $c_{d,k,t}(\varepsilon)$, with $c_{d,k,t}(0) = 0$, depending only on $d,k,t$ 
 and such that
 $$\mathrm{meas} \{ \lambda \in \order_w: \|\theta(\lambda)\| \leq \varepsilon \} \leq c_{d,k,t}(\varepsilon).$$
 \end{lemma}
 
 This follows from the fact that the space of $\theta$ is compact and
does not have $\theta = 0$ as a limit point; moreover, the assertion is true
``for each $\theta$ individually.''   We omit the easy formalization. 
%

\begin{lemma} \label{cdm}
Fix $d,k$. 
There is a continuous function $c_{k,d}(x)$ with $c_{k,d}(0) = 0$ which has the following property. 

Fix $\varepsilon > 0$.
 Let $\theta: F_w \rightarrow F_w^k$ be a nonconstant polynomial map of degree $d$,
and $x \in F_w$ so that $\|\theta(x)\| \leq \varepsilon$. Then there is a ball $B_x$ containing $x$ such that $\sup_{\lambda \in B_x} \|\theta(\lambda)\| \leq 1$ and 
$\meas(\lambda \in B_x: \|\theta(\lambda)\| \leq \varepsilon) \leq c_{k,d}(\varepsilon) \meas(B_x)$. 
\end{lemma}

\proof
Choose a maximal ball $B_x$ containing $x$ that satisfies the condition $\sup_{\lambda \in B_x}
\|\theta(\lambda)\| \leq 1$. We claim that there is a constant $c'$ depending only on $k,d$
so that $\sup_{\lambda \in B_x} \|\theta(\lambda)\| \geq c'$; this follows from
the interpolation formula
\begin{equation} \label{eq:interp}\theta(x) = \sum_{i=1}^{d+1} \theta(x_i) \frac{\prod_{j \neq i} (x-x_j)}{\prod_{j \neq i} (x_i - x_j)}\end{equation}

For simplicity let us assume that $d+1$ is smaller than the residue characteristic $q_w$ of $F_w$, the general case being similar. Suppose $B_x$ is the ball $\{\lambda: |\lambda - x|_w \leq q_w^{-K} \}$. 
We may choose $\{x_1, \dots, x_{d+1} \} \in B_x$
so that $|x_i - x_j|_w = q_w^{-K}$. 
On the other hand if $\lambda$ belongs to the enlarged ball $B' = \{\lambda: |\lambda - x|_w \leq q_w^{1-K}\}$
then $|\prod_{i \neq j} (\lambda - x_j) |_w \leq  q_w^{-d (K-1)}$. From this we see that 
$\sup_{\lambda \in B'} \|\theta\| \leq c' q_w^{d}$, which will contradict the maximality of $B'$ if $c'$ is too small. 

Now we apply the previous Lemma to the map $\theta$, rescaled
so that it is regarded as a map from $\order_w$ to $\order_w^k$.  \qed

\subsection{Convergence of limit measures.} \label{conditions}
Now let $\G$ be an algebraic group over a number field $F$, and $w$ a nonarchimedean place of $F$.
Set $G = \G(F_w)$, let $\Gamma$ be an arithmetic lattice of $G$, and let 
$\H$ be an algebraic subgroup of $\G$, $H = \H(F_w)$. 

We do not strive for generality and make the following assumptions, which are satisfied
in the context of the application in the text:
\begin{enumerate}
\item $\G$ is anisotropic over $F$, so that $\Gamma \backslash G$ is compact. 
(The general case is treated by a suitable variant of \cite[Theorem, \S 11.6]{MT}. )
\item The Lie algebra $\Lie(G)$ is simple as a Lie algebra over $F_w$. 
\item The subgroup $G^{+} \subset G$ generated by $F_w$-points of unipotent radicals
of parabolic $F_w$-subgroups coincides with $G$.  (If this fails, the statements must be slightly modified accordingly). 
\item $\H$ is semisimple, $H$ has property (T) and the subgroup $H^{+} \subset H$ generated by  $F_w$-points of unipotent radicals of parabolic $F_w$-subgroups coincides with $H$. (The last condition is absolutely essential, of course). 
\end{enumerate}
Let $\dot{\xi}_i \in \Gamma \backslash G$ be so that the orbits $\dot{\xi}_i H$ are closed; let $\mu_i$
be the $H$-invariant probability measure on $\dot{\xi}_i H$.

\begin{prop} \label{prop:converge}
Suppose $\xi_i \in G$ remain within a compact set and, for any subsequence of $i$,
the subgroups $\xi_i H \xi_i^{-1}$ generate\footnote{It will suffice that they generate a Zariski-dense subgroup of $G$, as will be clear from the proof. } $G$. Then any weak limit of the measures $\mu_i$ is the $G$-invariant probability measure on $\Gamma \backslash G$. 
\end{prop}

For any closed subgroups $H,L \subset G$, 
let $$X(H,L) = \{g \in G:  \mathrm{Ad}(g) \Lie(H)  \subset \Lie(L)\}.$$

We deduce Prop. \ref{prop:converge} from
\begin{lemma}\label{lem:converge} Notation being as in Prop. \ref{prop:converge},
assuming by passing to a subsequence that $\mu_i \rightarrow \mu_{\infty}$. 

Let $L$ be a proper subgroup of $G$ containing $H$, so that\footnote{On account of the assumption that $G^{+} = G$ and using the simplicity of $\Lie(G)$, any proper unbounded subgroup has lower dimension than $G$, as may be deduced from a theorem of Tits, see \cite{prasad}.} $\dim(L) < \dim(G)$. Let $\dot{\eta} \in \Gamma \backslash G$
so that $\dot{\eta} L$ is closed and supports an $L$-invariant probability measure.

Then there is a compact subset $X_0(H,L) \subset X(H,L)$ so that either:
\begin{enumerate} \item  \label{hyp1}
For infinitely many $i$, $\dot{\xi}_i H$ is contained in $\dot{\eta} L . X_0(H,L)$, or
\item \label{hyp2}  $\mu_{\infty}(\dot{\eta} L) = 0$. 
\end{enumerate}
\end{lemma}

Let us first deduce Prop. \ref{prop:converge} from Lem. \ref{lem:converge}. By Lem. \ref{ergodiclimit},
$\mu_{\infty}$ is an ergodic\footnote{In fact, Lemma \ref{ergodiclimit} was proved here only for $H$ a discrete group. Although Lemma \ref{ergodiclimit} is valid in general, as is shown in \cite{GW}, let us explicate how to obtain the desired conclusion in our context from this weaker form. In the present context, the fact that $H$ has compact orbits on $\Gamma \backslash G$
implies that $H$ admits a lattice $\Lambda \subset H$; then $\Lambda$ also has property (T), which is inherited by lattices. Each $\mu_i$ is $H$-ergodic and so also (by Howe-Moore) $\Lambda$-ergodic. Applying Lemma \ref{ergodiclimit} shows that $\mu_{\infty}$
is $\Lambda$-ergodic, so also $H$-ergodic.}
 $H$-invariant measure. By the measure classification theorem
of Ratner \cite{Ratner} and Margulis/Tomanov \cite{MT}, $\mu_{\infty}$ is {\em algebraic}: it is the $L$-invariant measure
supported on the closed subset $\dot{\eta} L \subset \Gamma \backslash G$, where $L \supset H$ is a closed subgroup and $L$ is the stabilizer of $\mu$ in $G$.  

It suffices to show that $L = G$. Suppose otherwise. 
Since the $\xi_i$ belong to a compact set and $\etadot L$ is compact, Lemma \ref{lem:converge} demonstrates that (after passing to a subsequence of $i$)
there is a compact subset $C \subset X(H,L)$, and a finite set $F \subset \Gamma$
such that $\xi_{i} \in F. \eta . C$ for all $i$. Passing to a further subsequence of $i$, we may assume that
there is a fixed $\gamma \in \Gamma$ so that $\xi_{i} \in \gamma \eta C$ for all $i$. 
Then $$\mathrm{Ad}(\xi_{i}) \Lie(H)  \subset \mathrm{Ad}(\gamma \eta) \Lie(L) $$
 for all $i$.  In particular, $\mathrm{Ad}(\xi_{i}) H $, in its adjoint action on $\Lie(G)$,  preserves $\mathrm{Ad} (\gamma \eta) \Lie(L)$. (The passage from $\Lie(H)$ to $H$ is effected using the fact that $H$ is generated by unipotent subgroups). The assumption on generation shows that $G$ preserves $\mathrm{Ad}(\gamma \eta) \Lie(L)$ also; since $\Lie(G)$ was assumed simple, this shows that $\Lie(L) = \Lie(G)$. This concludes the proof of Prop. \ref{prop:converge}. 
%



\proof (of Lem. \ref{lem:converge})




Let $\mathfrak{r}$ be a (vector space) complement to $\Lie(L)$ inside $\Lie(G)$ which
is stable by the conjugation action of $H$ (this is possible because, since $\H$ is semisimple, the adjoint action of $H$
on the Lie algebra is completely reducible).

   Let 
$\mathcal{\mathcal{B}}_1$ be an open compact neighbourhood of $0$ in $\mathfrak{r}$,
and $\mathcal{B}_r = \varpi_w^{r-1} \mathcal{B}_1$, for $r \geq 1$. 
We may assume that $\mathcal{B}_1$ is sufficiently small that the exponential
map\footnote{Which maps a neighbourhood of $0$ in the Lie algebra into $G$, equivariantly for the conjugation of $G$} is well-defined on $\mathcal{B}_1$, and moreover the map
$(\etadot L) \times \mathcal{B}_1 \rightarrow \Gamma \backslash G$
given by $(x, r) \mapsto x \exp(r)$ is a homeomorphism onto
an open neighbourhood $\mathcal{N}_1$ of $\etadot L$.
Define $\mathcal{N}_r$ to be the image of $\etadot L \times \mathcal{B}_r $
under this map.  Let $\pi: \mathcal{N}_1 \rightarrow \mathcal{B}_1$ be the natural projection
map, so that $\mathcal{N}_r = \pi^{-1} \mathcal{B}_r$. Set $X_0(H,L) = X(H,L) \cap \exp(\mathcal{B}_1)$.

Let $U \subset H$ be a one-parameter unipotent subgroup and $\theta:F_w \rightarrow U$
an isomorphism. 


By the  ergodicity of the $U$-action on $\xidot_i H$, 
for measure-generic points $y_i \in \xidot_i H$ 
the limit measure of the trajectory $y_i U$ 
is the measure $\mu_i$ for all $i$ (i.e., 
the $\theta(\order_w[-M])$-invariant probability measure on $y_i \theta(\order_w[-M])$
approaches the $H$-invariant probability measure on $\xidot_i H$, as $M \rightarrow \infty$). 
For such $y_i$, the closure $\overline{y_i U}$ coincides with $\xidot_i H$. 



Suppose $y_i$ is generic and belongs to $\mathcal{N}_1$. (If such does not exist, then $\mu_i(\mathcal{N}_1) = 0$ and we are done immediately). 
We may write $y_i = x_i \exp(r_i)$ for some $x_i \in \etadot L, r_i \in \mathcal{B}_1$.

So, for $\lambda \in F_w$, we have \begin{multline} \pi(y_i \theta(\lambda))) = 
\pi(x_i \exp(r_i) \theta(\lambda)) \\ = \pi(x_i \theta(\lambda)   \exp(\Ad \circ \theta(-\lambda) r_i)) = 
\Ad \circ \theta(-\lambda) r_i,\end{multline} so long as
$\Ad \circ \theta(-\lambda) r_i \subset \mathcal{B}_1$. 

Next we claim that either $\xidot_i H \subset \etadot L. X_0(H,L)$;  or the map $\lambda \mapsto \Ad \circ \theta(-\lambda) r_i$, which is visibly
polynomial from $F_w$ to $\mathfrak{r}$, is nonconstant for such $y_i$. Indeed, the closure $\overline{y_i U}$ is precisely
$\xidot H$, so, if $\lambda \mapsto \Ad \circ \theta(-\lambda) r_i$ were constant, we would have
in particular $$\xidot_i H \subset \etadot L .  \exp(r_i) = \etadot .  \exp(r_i) . (\exp(r_i)^{-1} L \exp(r_i)).$$  This implies that the Lie algebra of $H$ is contained
in the Lie algebra of $\exp(r_i)^{-1} L \exp(r_i)$. Therefore $\exp(r_i) \in X_0(H,L)$ and 
$\xidot_i H \subset \etadot L. X_0(H,L)$; so we are in the first case mentioned in the  Lemma. 


%
%

 Otherwise, set $Z_l = \{\lambda \in F_w: y_i \theta(\lambda) \in \mathcal{N}_l\}$, so that $F_w \supset
Z_1 \supset Z_2 \supset \dots$.  We note that the points $y_i \theta(\lambda)$ are generic
(in the sense above that their $U$-orbit is equidistributed w.r.t. $\mu_i$) for {\em all} $\lambda \in F_w$. 
Applying Lemma \ref{cdm} to the maps $\lambda \mapsto \pi(y_i \theta(\lambda))$, we see that, given $\varepsilon > 0$,
there exists $M$ big enough so that we can cover $Z_M$
by balls $B_j$ all contained in $Z_1$, and so that 
$\meas(B_j \cap Z_M)/\meas(B_j) \leq \varepsilon$ for each ball. 

It follows that, given any ball $Q \subset F_w$, there is a larger
ball $Q'$ such that $\meas(Q' \cap Z_M)/\meas(Q' ) \leq \varepsilon$.
(Either each ball $B_j$ corresponding to points in $Q \cap Z_M$
is contained in $Q$, or one such ball $B_{j_0}$ contains $Q$. In the former case, 
note that the family of maximal balls in the collection $\{B_j\}$ are disjoint and cover $Q \cap Z_M $;
take $Q' = Q$. In the latter case take $Q' = B_{j_0}$.) 
So the limit measure of the trajectory $y_i U$ assigns mass
$\leq \varepsilon$ to the neighbourhood $\mathcal{N}_M$.

Thus, if hypothesis (\ref{hyp1}) of the Lemma is not satisfied, we must have $\mu_i(\mathcal{N}_M) \leq \varepsilon$, for all $i$;
so the same is true for $\mu_{\infty}$ and so $\mu_{\infty}(\etadot L) = 0$, as required. \qed


\begin{thebibliography}{1}
\bibitem{Artin} M. Artin.
\newblock Geometric algebra. 
\newblock  Interscience Publishers, Inc., New York-London, 1957. 

\bibitem{Bhargava} M. Bhargava. 
\newblock Higher composition laws.
\newblock Princeton University PhD thesis, 2001. 

\bibitem{Bhargava2} M. Bhargava.
\newblock Higher composition laws. I. 
\newblock {\em Annals of Math.}, 159, 217--250, 2004. 

\bibitem{Bhargava3} M. Bhargava.
\newblock Higher composition laws. II. 
\newblock {\em Annals of Math.}, 159, 865--886, 2004. 

\bibitem{cassels} J. W. S. Cassels. 
\newblock Rational quadratic forms.
\newblock Academic Press, 1978. 

\bibitem{DM}
S. Dani and G. Margulis.
\newblock Orbits of unipotent flows and values of quadratic forms. 
\newblock {\em Advances in Soviet Mathematics}, vol 16, part 1. 
\newblock American Mathematical Society, 1993. 

\bibitem{Du}
W. Duke and R. Schulze-Pillot.
\newblock Representation of integers by positive ternary quadratic forms and equidistribution of lattice points on ellipsoids. 
\newblock {\em Invent. Math.} 99, 49--57, 1990. 


\bibitem{eichler}
M. Eichler.
\newblock 
Quadratische Formen und orthogonale Gruppen. 
\newblock Die Grundlehren der mathematischen Wissenschaften in Einzeldarstellungen mit besonderer BerŸcksichtigung der Anwendungsgebiete. Band LXIII.
Springer-Verlag, Berlin-Gšttingen-Heidelberg, 1952.

\bibitem{EL}
M. Einsiedler and E. Lindenstrauss. 
\newblock Diagonalizable flows on locally homogeneous spaces and number theory. 
\newblock ICM Proceedings 2006, to appear. 



\bibitem{ELMV}
M. Einsiedler, E. Lindenstrauss, Ph. Michel and A. Venkatesh. 
\newblock Distribution properties of compact orbits on homogeneous spaces, I, II \& III. 
\newblock In preparation. 



\bibitem{EMS}
A. Eskin, S. Mozes and N. Shah.
\newblock Unipotent flows and counting lattice points on homogeneous varieties. 
\newblock {\em Annals of Math.}, 143, 253--299, 1996. 


\bibitem{EO}
A. Eskin and H. Oh.
\newblock  Representations of integers by an invariant polynomial and unipotent flows.
\newblock Preprint, available {\tt http://www.its.caltech.edu/\~{}heeoh}. 

\bibitem{GG}
W. Gan and B. Gross. 
\newblock Commutative subgroups of certain non-associative rings. 
\newblock {\em Math. Ann.} 314, 265--283, 1999. 

\bibitem{GW}
Y. Glasner and B. Weiss.
\newblock Kazhdan's property T and the geometry of the collection of invariant measures. 
\newblock {\em Geom. Funct. Anal.}, 7, 917--935, 1997. 

\bibitem{Hsia}
J. Hsia.
\newblock Representations by spinor genera. 
\newblock {\em Pacific. J. Math.} 63, 147--152, 1976. 

\bibitem{hsia:2m+3}
J. Hsia, Y. Kitaoka and M. Kneser. 
\newblock Representations of positive definite quadratic forms.
\newblock {\em J. Reine Angew. Math.}, 301, 132--141, 1978. 

\bibitem{knes:jnt}
M. Kneser.
\newblock Composition of binary quadratic forms.
\newblock {\em J. Number Theory}, 15, 406--413, 1982. 

\bibitem{Kneser}
M. Kneser.
\newblock Darstellungsmasse indefiniter quadratischer Formen. 
\newblock {\em Math. Z}, 77, 188--194, 1961. 

\bibitem{linnik} 
Yu. Linnik.
\newblock Ergodic properties of algebraic number fields. 

\bibitem{MT}
G. Margulis and G. Tomanov.
\newblock Invariant measures for actions of unipotent groups over local fields on homogeneous spaces.
\newblock {\em Inventiones Math.}, 116, 347--392, 1994. 

\bibitem{MV}
\newblock Ph. Michel and A. Venkatesh.
\newblock Equidistribution of Gauss/Gross points to large moduli.
\newblock In preparation.    


\bibitem{Oh}
\newblock H. Oh.
\newblock Hardy-Littlewood system and representations of integers by an invariant polynomial. 
{\em Geom. Funct. Anal.} 14, 791--809, (2004). 

\bibitem{OMeara}
\newblock T. O'Meara.
\newblock Introduction to quadratic forms. 
\newblock 


\bibitem{Platonov}
V. P. Platonov.
\newblock The problem of strong approximation and the Kneser-Tits hypothesis for algebraic groups. (Russian)
\newblock {\em Izv. Akad. Nauk SSSR Ser. Mat.} 33, 1211--1219, 1969.


\bibitem{prasad}
G. Prasad. 
\newblock Elementary proof of a theorem of Bruhat-Tits-Rousseau and of a theorem of Tits. 
\newblock {\em Bull. Soc. Math.France}, 110, 197--202, 1982. 

\bibitem{Ratner}
M. Ratner.
\newblock Raghunathan's conjectures for Cartesian products of real and $p$-adic Lie groups.
\newblock {\em Duke. Math. J.}, 77, 275--382, 1995. 

\bibitem{shimura} G. Shimura. 
\newblock {\em Arithmetic and analytic theories of quadratic forms and Clifford groups.}
\newblock AMS, Providence, RI, 2004. 

\bibitem{schu:qfsurvey} R. Shulze-Pillot.
\newblock Representation by integral quadratic forms---a survey.  
\newblock {\em Algebraic and arithmetic theory of quadratic forms},  303--321, Contemp. Math., 344, 
\newblock Amer. Math. Soc., Providence, RI, 2004.
 

\bibitem{Tits} J. Tits.
\newblock Algebraic and abstract simple groups.
\newblock {\em Annals of Math.} (2), 80, 313--329, 1964.

\bibitem{weil} A. Weil.
\newblock Sur la th{\'e}orie des formes quadratiques.
\newblock {\em Colloq. ThŽorie des Groupes AlgŽbriques}, 9--22, Bruxelles, 1962.


\end{thebibliography}
\end{document}